# ON STATIONARITY OF LAGRANGIAN OBSERVATIONS OF PASSIVE TRACER VELOCITY IN A COMPRESSIBLE ENVIRONMENT[1]

By Tomasz Komorowski and Grzegorz Krupa

*University of Maria Curie-Sklodowska and Catholic University of Lubin*

We study the transport of a passive tracer particle in a steady strongly mixing flow with a nonzero mean velocity. We show that there exists a probability measure under which the particle Lagrangian velocity process is stationary. This measure is absolutely continuous with respect to the underlying probability measure for the Eulerian flow.

**1. Introduction.** One of the simplest models of the passive tracer motion in a turbulent flow is given by the Itô stochastic equation

$$\begin{aligned}
d\mathbf{x}(t) &= \mathbf{u}(\mathbf{x}(t))\,dt + \sqrt{\kappa}\,d\mathbf{w}(t), \qquad t \geq 0, \\
\mathbf{x}(0) &= \mathbf{0}.
\end{aligned} \tag{1.1}$$

Here $\mathbf{u} = (u_1, \ldots, u_d) : \mathbb{R}^d \times \Omega \to \mathbb{R}^d$ is the *Eulerian velocity field* of the flow. It is assumed to be a stationary, $d$-dimensional random vector field given over a certain probability space $(\Omega, \mathcal{V}, \mathbb{P})$, and $\mathbf{w}(\cdot)$ is a standard $d$-dimensional Brownian motion defined over another probability space $(\Sigma, \mathcal{A}, Q)$. Parameter $\kappa$, called the molecular diffusivity of the medium, is assumed to be nonnegative. The resulting process $\mathbf{x}(\cdot)$ is considered over the product probability space $(\Omega \times \Sigma, \mathcal{V} \otimes \mathcal{A}, \mathbb{P} \otimes Q)$. A question that generates considerable interest in statistical hydrodynamics is to provide the description of the long-time, large-scale asymptotics of $\mathbf{x}(\cdot)$. Possible types of the trajectory behavior that may occur in the limit include Newtonian motions, diffusions, fractional diffusions and possibly Lévy flights; see [2, 6, 7, 18].

Received March 2002; revised September 2003.
[1]Supported by the Polish Committee for Scientific Research (KBN) Grant 2PO3A03123.

An important insight in understanding the asymptotic behavior of solutions to (1.1) can be gained if one is able to establish the existence of a probability measure $\mu$ defined over $\mathcal{V} \otimes \mathcal{A}$ under which the *Lagrangian process*, that is, $\mathbf{u}(\mathbf{x}(t))$, $t \geq 0$, is stationary and ergodic. The above process corresponds to the observations of the velocity from the vantage point of the moving trajectory. If such a measure exists, then one can conclude that

$$(1.2) \qquad \mathbf{v}_* := \lim_{t \uparrow +\infty} \frac{\mathbf{x}(t)}{t} = \int \mathbf{u}(\mathbf{0}) \, d\mu$$

exists $\mu$-a.s. and $\mathbf{v}_*$ is deterministic. $\mathbf{v}_*$ is sometimes called the *Stokes drift* of the medium. If, in addition to stationarity and ergodicity, $\mu$ is absolutely continuous w.r.t. the product measure $\mathbb{P} \otimes Q$ and the limit in (1.2) holds $\mathbb{P} \otimes Q$-a.s., we shall call $\mu$ *a regular, invariant measure* for the Lagrangian process.

In the present paper we consider the case of strongly mixing, steady (i.e., time-independent) velocity fields. The main result we set out to prove can be stated informally as follows; see Theorem 2.2 for the precise statement.

THEOREM 1.1. *Suppose that the molecular diffusivity $\kappa$ is strictly positive, and that the velocity field $\mathbf{u}$ is stationary with the mean that is larger than the amplitude of its fluctuations [see condition (A)] and decorrelates at finite distances [condition (DR)]. Then, assuming also some topological and measure-theoretic regularity properties of the field [condition (R)], there exists a regular invariant measure $\mu$ for the Lagrangian process $\mathbf{u}(\mathbf{x}(t))$, $t \geq 0$.*

The standard results, for example, those obtained in the framework of the homogenization theory (see [13]), concern the drifts that are either gradients of stationary scalar potentials (i.e., $\mathbf{u} = \nabla_{\mathbf{x}}\phi$, where $\phi$ is a certain stationary scalar field), or are incompressible (i.e., $\nabla_{\mathbf{x}} \cdot \mathbf{u} := \sum_{i=1}^{d} \partial_{x_i} u_i \equiv 0$). The gradient case corresponds to the motion of a tracer particle in a medium (e.g., gas) that remains in an equilibrium, while the incompressible fields describe turbulent flows of fluids. In both of these cases, regular invariant measures $\mu$ can be identified explicitly. In the gradient case $\mu$ is given by $\mathbb{P}_\phi \otimes Q$, where $\mathbb{P}_\phi$ the is Gibbs equilibrium measure relative to the potential $\phi/\kappa$, while in the incompressible case the invariant measure is actually equal to $\mathbb{P} \otimes Q$; see [[12]–[14]].

In many interesting situations, however, the motion of a tracer takes place in a compressible environment that is far from being in equilibrium, for example, floating of a particle on a fluid surface; see [5]. Due to the infinite-dimensional character of the problem, the existence of a regular invariant measure is, in general, hard to prove and very few results concerning the problem are known. For a review of the existing literature on the subject, a reader is advised to consult [19]. It is generally believed, however, that strong mixing properties of the Eulerian flow should guarantee the existence of a regular invariant measure for Lagrangian observations [5, 19]. Recently, a number of rigorous results substantiating that point of view have been obtained for nonsteady (time-dependent) flows; see [[8]–[10]]. A generic situation considered in those papers concerns fields that have strong temporal decorrelation properties; for example, in [8, 9] the Eulerian velocity field is of finite time dependence range, while in [10] it is Gaussian, Markovian and sufficiently strongly mixing in the temporal variable. Under any of these assumptions, it can be shown that there exists then a regular (w.r.t. $\mathbb{P} \otimes Q$) invariant measure $\mu$, provided the molecular diffusivity $\kappa$ is positive. It is worthwhile to point out that under some additional assumptions on the spatial structure of the velocity field (i.e., the compact support of the spatial power-energy spectrum), the results of [10] deal also with the case of the vanishing molecular diffusivity $\kappa = 0$.

Let us discuss briefly the principal ideas of the proof of Theorem 1.1. First, we use the factorization property of the $\sigma$-algebra corresponding to the Eulerian velocity field. More specifically, looking in the direction $\hat{\mathbf{v}}$ pointed by the mean drift of the flow at any given time instant, say $t = 0$, we can decompose the future history of the velocity field into the part that is determined by the tracer particle history up to $t = 0$ and, independent of it, the renewal part. This decomposition forms the base for the definition of the so-called *transport operator*; see Section 4. Informally speaking, it describes the change of the statistics of the field, as observed from the moving particle, within the time span $\tau$ needed for the particle to travel from the initial position at $t = 0$ the spatial distance required for the complete renewal of the Eulerian velocity. In addition, after this time the particle does not revisit the half-space containing the initial position of the trajectory and bounded by the hyperplane orthogonal to the drift passing by the point which is unit to the left in the direction $\hat{\mathbf{v}}$ from the position of the particle at $\tau$. Because of this property, we call $\tau$ the *nonretraction time*; see Section 3 for its precise definition. It is obviously a non-Markovian random time. The definition of $\tau$ is modelled on the notion of the nonretraction times, introduced by Sznitman and Zerner in [17], in the case of random walks on a random integer lattice with independent sites. As we show in Section 4, see (4.3), the transport operator acts on a certain space of density functions with respect to $\mathbb{P}$. An important property of this operator is the fact that it admits an invariant density; see Proposition 4.7. This density is subsequently used, see formula (5.33), to define the invariant measure $\mu$; see Theorem 1.1.

A result, that corresponds to our main theorem has been proved for the nearest-neighbor random walks on integer lattice $\mathbb{Z}^d$ with i.i.d. transition probabilities having a uniform local drift property (the so-called nonnestling condition) in [3].

**2. Notation and formulation of the main result.** To simplify the notation we assume, throughout the remainder of the paper, that $\kappa = 1$ in (1.1). The proof can be trivially generalized to the case of an arbitrary positive molecular diffusivity. The case $\kappa = 0$ is substantially different and we know of no results concerning the steady fields in that situation. (As mentioned earlier, some results concerning time-dependent, Gaussian, Markovian fields can be found in [10].)

For any $L > 0$, we denote $\mathfrak{X}_L := C([0,L];\mathbb{R}^d)$ and $\mathfrak{X} := C([0,+\infty);\mathbb{R}^d)$. These spaces are equipped with the standard topology of uniform convergence on compact sets. For any $t \geq 0$, we denote by $\Pi(t) \colon \mathfrak{X} \to \mathbb{R}^d$ the canonical projection $\Pi(t)(\pi) := \pi(t)$, $\pi \in \mathfrak{X}$. Let $\mathcal{M}_t := \sigma\{\Pi(s) : s \leq t\}$, $t \geq 0$, be the canonical filtration on $\mathfrak{X}$. We let $\mathcal{M} := \bigvee_{t \geq 0} \mathcal{M}_t$. By $\mathcal{P}$ and $\mathcal{P}_L$ we denote the spaces of all Borel probability measures on $\mathfrak{X}$ and $\mathfrak{X}_L$, respectively. By $\mathbb{W}$ and $\mathbb{W}_L$ we denote the standard Wiener measure on $(\mathfrak{X}, \mathcal{M})$ and its restriction to $\mathcal{M}_L$, respectively. For any $h \geq 0$, we define the shift operator $\theta_h \colon \mathfrak{X} \to \mathfrak{X}$ given by $\theta_h(\pi)(t) := \pi(t+h)$ for all $t \geq 0$, $\pi \in \mathfrak{X}$.

Let $(\Omega, d)$ be a Polish space. We denote by $\mathcal{B}(\Omega)$ the $\sigma$-algebra of Borel subsets of $\Omega$. We suppose that $\mathbb{P}$ is a Borel probability measure and $\mathbb{E}[\cdot]$ denotes the corresponding expectation. Let $\mathcal{N}$ be the $\sigma$-ring of $\mathbb{P}$-null sets of $\overline{\mathcal{B}(\Omega)}$, the $\mathbb{P}$-completion of $\mathcal{B}(\Omega)$. Unless otherwise stated, we will assume that any sub-$\sigma$-algebra of $\overline{\mathcal{B}(\Omega)}$ under consideration contains $\mathcal{N}$. For brevity we write $L^p := L^p(\mathcal{T}_1)$, $p \in [1, +\infty]$, where $\mathcal{T}_1 := (\Omega, \overline{\mathcal{B}(\Omega)}, \mathbb{P})$. We assume further the property of spatial homogeneity of measure $\mathbb{P}$. The above means that there exists a group of transformations $T_{\mathbf{x}}, \mathbf{x} \in \mathbb{R}^d$, acting on $\Omega$ such that, for any $\mathbf{x} \in \mathbb{R}^d$, $A \in \overline{\mathcal{B}(\Omega)}$, we have $T_{\mathbf{x}}(A) \in \overline{\mathcal{B}(\Omega)}$ and $\mathbb{P}(T_{\mathbf{x}}(A)) = \mathbb{P}(A)$.

We assume that $\mathbf{u} \colon \Omega \to \mathbb{R}^d$ is a random vector over $\mathcal{T}_1$ satisfying

(A) $|\mathbf{v}| > \|\widetilde{\mathbf{u}}\|_{L^\infty}$, where $\mathbf{v} := \mathbb{E}\mathbf{u}$ and $\widetilde{\mathbf{u}} = \mathbf{u} - \mathbf{v}$.

The spatially homogeneous Eulerian velocity field is defined as $\mathbf{u}(\mathbf{x}; \omega) := \mathbf{u}(T_{\mathbf{x}}(\omega))$. Assumption (A) guarantees that the mean drift dominates its fluctuations and therefore there exists $\delta > 0$ such that

$$(2.1) \qquad \mathbf{u}(\mathbf{x}) \cdot \hat{\mathbf{v}} \geq \delta > 0, \qquad \mathbb{P}\text{-a.s.}$$

for all $\mathbf{x} \in \mathbb{R}^d$. Here $\hat{\mathbf{v}} := \mathbf{v}/|\mathbf{v}|$. We shall also assume that $1 \wedge (|\mathbf{v}|/4) \geq \delta > 0$.

For any $R > 0$, we denote by $\mathcal{F}_R^i$, $\mathcal{F}_R^e$ the $\sigma$-algebras generated by $\mathbf{u}(\mathbf{x})$, $|\mathbf{x}| \leq R$, and $\mathbf{u}(\mathbf{x})$, $|\mathbf{x}| \geq R$, correspondingly. We assume that

(DR) (*finite dependence range*) there exists $r_0 > 0$ such that, for any $r > 0$, the $\sigma$-algebras $\mathcal{F}_r^i$ and $\mathcal{F}_{r+r_0}^e$ are independent.

Finally, we suppose that the field possesses certain regularity both in the topological and the measure-theoretic sense. Namely, we assume that

(R) for any $\omega \in \Omega$, the field $\mathbf{u}(\cdot;\omega)$ is of $C^1$ class of regularity and there exists a deterministic constant $U > 0$ such that $\|\tilde{\mathbf{u}}(\cdot;\omega)\|_{W^{1,\infty}(\mathbb{R}^d)} \leq U$. The norm taken here is the usual one corresponding to the classical Sobolev space $W^{1,\infty}(\mathbb{R}^d)$.

In addition, we suppose that, for any $N \geq 1$ and $\mathbf{x}_1, \ldots, \mathbf{x}_N \in \mathbb{R}^d$ such that $\mathbf{x}_i \neq \mathbf{x}_j$ for $i \neq j \in \{1, \ldots, N\}$, the probability distribution of $(\mathbf{u}(\mathbf{x}_1), \ldots, \mathbf{u}(\mathbf{x}_N))$ in the space $(\mathbb{R}^d)^N$ is absolutely continuous with respect to the $N \cdot d$-dimensional Lebesgue measure.

REMARK 2.1. Let us briefly discuss some important, from our point of view, consequences of assumption (R). Let $t \in \mathbb{R}$ and $\mathcal{V}_t$ be the sub-$\sigma$-algebra of $\overline{\mathcal{B}(\Omega)}$ generated by $\mathbf{u}(\mathbf{x})$, $\mathbf{x} \cdot \hat{\mathbf{v}} \leq t$. We note that obviously

$$T_{\mathbf{x}}(\mathcal{V}_t) \subseteq \mathcal{V}_{t+\mathbf{x}\cdot\hat{\mathbf{v}}} \qquad \forall\, (t, \mathbf{x}) \in \mathbb{R} \times \mathbb{R}^d.$$

Thanks to the assumption (R) (see page 66 of [15], or the Appendix of [9]), the filtration $(\mathcal{V}_t)_{t \geq 0}$ admits a *factorization* with respect to $\mathcal{V}_0$, that is, for any $t \geq 0$, there exists a $\sigma$-algebra $\mathcal{R}^t$ such that $\mathcal{V}_0$ and $\mathcal{R}^t$ are $\mathbb{P}$-independent and $\mathcal{V}_t$ is generated by $\mathcal{V}_0$ and $\mathcal{R}^t$. Let $\mathcal{R} := \bigvee_{t \geq 0} \mathcal{R}^t$.

Note that $(\mathcal{R}^t)_{t \geq 0}$ form a filtration of $\sigma$-algebras. Indeed, any random variable $H(\cdot)$ that is $\mathcal{R}^t$-measurable is $\mathcal{V}_u$-measurable for any $u \geq t$ and one can find (from the factorization property) a random variable $G(\cdot, \cdot)$ that is $\mathcal{V}_0 \otimes \mathcal{R}^u$-measurable and $H(\omega) = G(\omega, \omega)$. From the fact that $H$ is independent of $\mathcal{V}_0$, we immediately conclude that $H(\cdot) = \int G(\omega', \cdot)\mathbb{P}(d\omega')$, $\mathbb{P}$-a.s., thus $H$ is $\mathcal{R}^u$-measurable.

The previous argument also shows that, thanks to condition (DR), any random vector $\mathbf{u}(\mathbf{x})$, with $t := \mathbf{x} \cdot \hat{\mathbf{v}} \geq r_0$ and $r_0$ as in (DR), is $\mathcal{R}^t$-measurable.

Let $Q_{\mathbf{x},\omega} \in \mathcal{P}$ be the law of the solution to (1.1) for a fixed realization of $\omega \in \Omega$ and subject to the initial condition $\mathbf{x}(0) = \mathbf{x}$. We denote $\mathcal{T}_{\mathbf{x},\omega} := (\mathfrak{X}, \mathcal{M}, Q_{\mathbf{x},\omega})$ and by $\mathbf{M}_{\mathbf{x},\omega}$ the respective mathematical expectation. In the particular case when $\mathbf{x} = \mathbf{0}$, we shall suppress the subscript $\mathbf{x}$.

The process

(2.2) $$\mathbf{w}_\omega(t; \pi) := \pi(t) - \int_0^t \mathbf{u}(\pi(s); \omega)\, ds, \qquad t \geq 0,$$

is a $d$-dimensional standard Brownian motion starting at $\mathbf{x}$ over $\mathcal{T}_{\mathbf{x},\omega}$ for any $\omega$. We denote by $p^\omega : \mathbb{R}_+ \times \mathbb{R}^d \times \mathbb{R}^d \to \mathbb{R}_+$ the transition probability densities of the diffusion given by (1.1).

Define a measure $P_0$ on $(\Omega \times \mathfrak{X}, \overline{\mathcal{B}(\Omega)} \otimes \mathcal{M})$ as the semiproduct

$$P_0(A \times B) := \int_A Q_\omega(B)\mathbb{P}(d\omega) \qquad \forall\, A \in \overline{\mathcal{B}(\Omega)},\ B \in \mathcal{M},$$

and a stochastic process

$$V(t;\omega,\pi) := \mathbf{u}(\pi(t);\omega), \qquad t \geq 0, \tag{2.3}$$

over $(\Omega \times \mathfrak{X}, \overline{\mathcal{B}(\Omega)} \otimes \mathcal{M}, P_0)$.

Theorem 1.1 can be stated more precisely in the following way.

THEOREM 2.2. *Suppose that* $\mathbf{u}$ *is a velocity field satisfying the assumptions* (A), (DR) *and* (R). *Then, there exists a probability measure* $\mu$ *on* $(\Omega \times \mathfrak{X}, \overline{\mathcal{B}(\Omega)} \otimes \mathcal{M})$ *that is absolutely continuous w.r.t.* $P_0$ *and such that the process* $V(\cdot)$ *is stationary and ergodic w.r.t.* $\mu$. *In addition, the law of large numbers holds w.r.t.* $P_0$, *that is,*

$$\lim_{t \uparrow +\infty} \frac{\pi(t)}{t} = \int \mathbf{u}(0) \, d\mu, \qquad P_0\text{-a.s.} \tag{2.4}$$

In Theorem 2.2, ergodicity of the relevant measure is understood as the absence of shift-invariant nontrivial sets. More precisely, any Borelian subset $A \subseteq \mathfrak{X}$ such that

$$\int |\mathbf{1}_{\theta_h(A)}(V(\cdot)) - \mathbf{1}_A(V(\cdot))| \, d\mu = 0 \qquad \text{for all } h \geq 0,$$

must be $\mu$-trivial, that is

$$\mu[(\omega,\pi) : V(\cdot;\omega,\pi) \in A] = 0 \text{ or } 1. \tag{2.5}$$

**3. Nonretraction times.** In this section we introduce a family of random variables that, for reasons which become obvious later on, we shall call *nonretraction times*. They are not stopping times and describe subsequent times after which no retraction of the diffusion can occur in the direction pointed by the mean velocity. This notion is based on a discrete analogue introduced for random walks on a random lattice in [17]. Since the results contained in this section are modifications of the corresponding results of [17], we postpone their proofs to Appendix A.

For any $\pi \in \mathfrak{X}$, $l \in [0, +\infty)$, we let

$$D(l;\pi) := \min[t \geq 0 : \hat{\mathbf{v}} \cdot \pi(t) \leq -1 + l]. \tag{3.1}$$

For brevity sake we write $D := D(\hat{\mathbf{v}} \cdot \pi(0))$. Let

$$U_u(\pi) := \min[t \geq 0 : \hat{\mathbf{v}} \cdot \pi(t) \geq u], \tag{3.2}$$

$$\tilde{U}_u(\pi) := \min[t \geq 0 : \hat{\mathbf{v}} \cdot \pi(t) \leq u]$$

and

$$M_*(\pi) := \sup[\hat{\mathbf{v}} \cdot (\pi(t) - \pi(0)) : 0 \leq t \leq D(\pi)]. \tag{3.3}$$

The last random variable is defined for those $\pi$, for which $D(\pi) < +\infty$.

For any $t \geq 0$, we define also

$$(3.4) \qquad A(t) := \left[\pi \colon \inf_{s \in [0,t]} (\pi(s) \cdot \hat{\mathbf{v}} - \pi(0) \cdot \hat{\mathbf{v}}) \geq -1\right].$$

We introduce the sequence of $(\mathcal{M}_t)$-stopping times $(S_k)_{k \geq 0}$, $(R_k)_{k \geq 0}$ and the sequence of maxima $(M_k)_{k \geq 0}$ letting

$$S_0 := 0, \qquad R_0 := 0, \qquad M_0 := \hat{\mathbf{v}} \cdot \pi(0),$$
$$(3.5) \qquad S_1 := U_{M_0 + r_0 + 1} \leq +\infty, \qquad R_1 := D \circ \theta_{S_1} + S_1 \leq +\infty,$$
$$M_1 := \max[\hat{\mathbf{v}} \cdot \pi(t), 0 \leq t \leq R_1] \leq +\infty,$$

where $r_0 > 0$ is as in (DR).

By induction we set, for any $k \geq 1$,

$$(3.6) \qquad \begin{aligned} S_{k+1} &:= U_{M_k + r_0 + 1}, \qquad R_{k+1} := D \circ \theta_{S_{k+1}} + S_{k+1}, \\ M_{k+1} &:= \max[\hat{\mathbf{v}} \cdot \pi(t), 0 \leq t \leq R_{k+1}]. \end{aligned}$$

The following lemmas hold.

LEMMA 3.1. *There exist deterministic constants $\gamma, \gamma_1, \gamma_2 > 0$ such that, for each $\mathbf{x} \in \mathbb{R}^d$, we have*

$$(3.7) \qquad Q_{\mathbf{x},\omega}[D = +\infty] \geq \gamma, \qquad \mathbb{P}\text{-}a.s.,$$

$$(3.8) \qquad Q_{\mathbf{x},\omega}[\tilde{U}_{\mathbf{x} \cdot \hat{\mathbf{v}} - M} < U_{\mathbf{x} \cdot \hat{\mathbf{v}} + M}] \leq \exp\{-\gamma_1 M\} \qquad \forall M > 0, \ \mathbb{P}\text{-}a.s.,$$

$$(3.9) \qquad \mathbf{M}_{\mathbf{x},\omega}[M_*, D < +\infty] \leq \gamma_2, \qquad \mathbb{P}\text{-}a.s.$$

LEMMA 3.2. *With the notation of Lemma 3.1, for each $\mathbf{x} \in \mathbb{R}^d$, we have*

$$(3.10) \qquad \limsup_{m \uparrow +\infty} \mathbf{M}_{\mathbf{x},\omega}\left[\frac{U_m}{m}\right] \leq \frac{1}{\delta}, \qquad \mathbb{P}\text{-}a.s.$$

*and*

$$(3.11) \qquad Q_{\mathbf{x},\omega}[R_k < +\infty] \leq (1-\gamma)^k \qquad \forall k \geq 1, \ \mathbb{P}\text{-}a.s.$$

Let $K := \inf[k \geq 1 \colon R_k = +\infty]$, or $K = +\infty$ if the set of which we take the infimum is empty.

COROLLARY 3.3. *For each $\mathbf{x} \in \mathbb{R}^d$, we have:*

(i) $Q_{\mathbf{x},\omega}[K < +\infty] = 1$, $\mathbb{P}$-*a.s. and*
(ii) $Q_{\mathbf{x},\omega}[S_K < +\infty] = 1$, $\mathbb{P}$-*a.s.*

We define the first nonretraction time $\tau_1 := S_K < +\infty$, $P_0$-a.s. The subsequent times of nonretraction $\tau_n$, $n \geq 2$, are defined by induction using the relation

$$\tau_{n+1} := \tau_n + \tau_1 \circ \theta_{\tau_n} \qquad \text{for } n \geq 1. \tag{3.12}$$

Note that the random variables $\tau_n$ need not be $(\mathcal{M}_t)$-stopping times.

**4. The transport operator and its properties.** For any $a, b \in \mathbb{R} \cup \{-\infty, +\infty\}$, $a \leq b$, we let $\mathcal{V}_{a,b}$ be the $\sigma$-algebra generated by $\mathbf{u}(\mathbf{x})$, where $a < \mathbf{x} \cdot \hat{\mathbf{v}} < b$. Consistently with Remark 2.1, we write $\mathcal{V}_a$ for $\mathcal{V}_{-\infty,a}$. Let $\mathcal{T}_2 := (\Omega, \mathcal{V}_0, \mathbb{P})$,

$$\mathbb{P}_D(d\omega) := \frac{Q_\omega[D = +\infty]}{P_0[D = +\infty]} \mathbb{P}(d\omega),$$

$$P_D(d\omega, d\pi) := \frac{\mathbf{1}_{[D(\pi) = +\infty]}}{P_0[D = +\infty]} P_0(d\omega, d\pi),$$

and $\mathcal{T}_D := (\Omega, \mathcal{V}_0, \mathbb{P}_D)$. Note that in light of (3.7), $\mathbb{P}_D$ is equivalent with $\mathbb{P}$. Also, for any probability triple $\mathcal{T}$ the symbol $\mathcal{D}(\mathcal{T})$ denotes the set of all probability densities w.r.t. the relevant probability measure, that is, non-negative elements of $L^1(\mathcal{T})$ whose integral equals 1.

In this section we introduce a certain linear operator $\mathfrak{Q} : L^1(\mathcal{T}_D) \to L^1(\mathcal{T}_D)$ that preserves $\mathcal{D}(\mathcal{T}_D)$. It is conjugate, in the sense of (4.1), with the spatial shift by $\pi(\tau_1)$, that is,

$$\int \mathbf{M}_\omega[G(T_{\pi(\tau_1)}(\omega)), D = +\infty] F(\omega) \mathbb{P}(d\omega)$$

$$= \int G(\omega) \mathfrak{Q} F(\omega) Q_\omega[D = +\infty] \mathbb{P}(d\omega) \tag{4.1}$$

(recall that $Q_\omega := Q_{\mathbf{0},\omega}$, $\mathbf{M}_\omega := \mathbf{M}_{\mathbf{0},\omega}$) for any $F$ and $G$ that are correspondingly $\mathcal{V}_0$ and $\mathcal{V}_{0,+\infty}$-measurable; see Proposition 4.6. We will call $\mathfrak{Q}$ a *transport operator*.

4.1. *Some consequences of the factoring property.* Let $\mathcal{T}_3 := (\Omega, \mathcal{R}, \mathbb{P})$ and let $\mathcal{T}_2 \otimes \mathcal{T}_3 := (\Omega \times \Omega, \mathcal{V}_0 \otimes \mathcal{R}, \mathbb{P} \otimes \mathbb{P})$. Condition (R) implies (see Appendix B) the existence of an isometric isomorphism $\mathcal{Z} : L^p(\mathcal{T}_1) \to L^p(\mathcal{T}_2 \otimes \mathcal{T}_3)$, $p \in [1, \infty]$ such that:

(Z1) $\mathcal{Z}F \geq 0$ for $F \geq 0$ and $\mathcal{Z}\mathbf{1} = \mathbf{1}$,
(Z2) for any $F_1, \ldots, F_N \in L^p(\mathcal{T}_1)$ and $\Phi : \mathbb{R}^N \to \mathbb{R}$ bounded and continuous, we have $\mathcal{Z}(\mathbf{\Phi}(F_1, \ldots, F_N)) = \mathbf{\Phi}(\mathcal{Z}F_1, \ldots, \mathcal{Z}F_N)$,
(Z3) $\mathcal{Z}F(\omega, \omega') = F(\omega)$ for all $F \in L^p(\mathcal{T}_2)$, $\mathcal{Z}G(\omega, \omega') = G(\omega')$ for all $G \in L^p(\mathcal{T}_3)$,
(Z4) $\mathcal{Z}F$ is $\mathcal{V}_0 \otimes \mathcal{R}^t$-measurable if $F$ is $\mathcal{V}_t$-measurable, for any $t \geq 0$.

REMARK 4.1. From condition (Z2) we conclude immediately the following:

(Z5) $\mathcal{Z}F$ has the same law as $F$ for all $F \in L^1(\mathcal{T}_1)$.
(Z6) Suppose that $F_1, F_2 \in L^\infty(\mathcal{T}_1)$. Then $\mathcal{Z}(F_1 F_2) = \mathcal{Z}(F_1)\mathcal{Z}(F_2)$.

REMARK 4.2. Denote $\mathbf{U}(\mathbf{x}) := \mathcal{Z}(\mathbf{u}(\mathbf{x})) \in L^\infty(\mathcal{T}_2 \otimes \mathcal{T}_3)$ for any fixed $\mathbf{x} \in \mathbb{R}^d$. One can find a modification of $\mathbf{U}$ defined over $\mathcal{T}_2 \otimes \mathcal{T}_3$ that is of $C^1$-class of regularity $\mathbb{P} \otimes \mathbb{P}$-a.s. and such that $\|\tilde{\mathbf{U}}(\cdot; \omega, \omega')\|_{W^{1,\infty}(\mathbb{R}^d)} \leq U$ [here $U$ is as in the statement of condition (R)] for $\mathbb{P} \otimes \mathbb{P}$-a.s. $(\omega, \omega')$, where $\tilde{\mathbf{U}}(\cdot; \omega, \omega') := \mathbf{U}(\cdot; \omega, \omega') - \mathbf{v}$.

4.2. *The definition of operator* $\mathfrak{Q}$. We start with some auxiliary notation. Let $\mathcal{T}_\mathbb{W} := (\mathfrak{X}, \mathcal{M}, \mathbb{W})$, where, as we recall, $\mathbb{W}$ is the standard Wiener measure. By $\mathbb{F}_k$ we denote the law in $\mathbb{R}^d$ of random vector $\pi(S_k)$ over $\mathcal{T}_\mathbb{W}$. Let

$$\nu_L(\pi; \omega) := \exp\left\{ \int_0^L \mathbf{u}(\pi(s); \omega) \, d\pi(s) - \tfrac{1}{2} \int_0^L |\mathbf{u}(\pi(s); \omega)|^2 \, ds \right\}$$

be the Radon–Nikodym derivative $\frac{dQ_{\omega,L}}{d\mathbb{W}_L}$. Here $Q_{\omega,L}$, $\mathbb{W}_L$ is the restriction of $Q_\omega$, $\mathbb{W}$ to $\mathcal{M}_L$ for a given $L > 0$. $\int_0^t \mathbf{u}(\pi(s); \omega) \, d\pi(s)$, $t \geq 0$, is the stochastic integral with respect to the Wiener process $\pi(\cdot)$ over the probability space $\mathcal{T}_\mathbb{W}$, (see [16], page 99). Using the properties (Z1), (Z2) and (Z6) of the operator $\mathcal{Z}$, one can prove that

(4.2) $$\mathcal{Z}(\nu_L(\pi))(\omega, \omega') = \bar{\nu}_L(\pi; \omega, \omega'),$$

where

$$\bar{\nu}_L(\pi; \omega, \omega') := \exp\left\{ \int_0^L \mathbf{U}(\pi(s); \omega, \omega') \, d\pi(s) - \tfrac{1}{2} \int_0^L |\mathbf{U}(\pi(s); \omega, \omega')|^2 \, ds \right\}$$
$$\forall L > 0.$$

For $\mathbf{x} \in \mathbb{R}^d$, we define $\mathbb{W}_{k,L,\mathbf{x}} \in \mathcal{P}_L$ and $\mathbf{M}_{k,L,\mathbf{x}}$, the regular conditional probabilities obtained by conditioning of $\mathbb{W}_L$ on the event $\{\pi(S_k) = \mathbf{x}, S_k \in [L-1, L)\}$ and the respective expectations.

The linear operator $\mathfrak{Q}$ satisfying (4.1) will be defined as follows. For any $F$ that is bounded and $\mathcal{V}_0$-measurable, define

(4.3) $$\mathfrak{Q}F(\omega') := \int \mathcal{K}(\omega, \omega') F(\omega) \mathbb{P}(d\omega),$$

where

(4.4) $$\mathcal{K}(\omega, \omega') := \sum_{k,L=1}^{+\infty} \int_{\mathbb{R}^d} \mathbf{M}_{k,L,\mathbf{x}}(\omega, T_{-\mathbf{x}}\omega') \mathbb{F}_k(d\mathbf{x})$$

and

(4.5) $\quad \mathcal{M}_{k,L,\mathbf{x}}(\omega,\omega') := \mathbf{M}_{k,L,\mathbf{x}}[\bar{\nu}_{S_k}(\pi;\omega,\omega'), A(S_k), L-1 \leq S_k < L].$

Let
$$t_k(\pi) := \hat{\mathbf{v}} \cdot \pi(S_k).$$

Note that $\bar{\nu}_{S_k}$ is $\mathcal{V}_0 \otimes \mathcal{R}^{t_k(\pi)}$-measurable for $\mathbb{W}_L$-a.s. $\pi$. Hence $\mathcal{M}_{k,L,\mathbf{x}}(\cdot, T_{-\mathbf{x}}\cdot)$ is $\mathcal{V}_0 \otimes \mathcal{V}_0$-measurable for $\mathbb{F}_k$-a.s. $\mathbf{x}$.

REMARK 4.3. The definition of the operator $\mathfrak{Q}$ given by (4.3) and (4.4) may seem to look a bit technical at the moment. To motivate it we remark here that $\mathfrak{Q}$ is constructed in such a way that the property expressed by (4.1) holds; see part (i) of Proposition 4.5 and Proposition 4.6. This property enables us to reduce the question of the existence of an absolutely continuous modification of measure $P_0$ under which the sequence $(\tau_{k+1} - \tau_k, \pi(\tau_{k+1}) - \pi(\tau_k))$, $k \geq 1$, is stationary to the problem of the existence of an invariant density for $\mathfrak{Q}$; see Theorems 4.7 and 5.1. This result together with integrability of the moments of $\tau_1$ and $|\pi(\tau_1)|$ (see Proposition 5.2) allow us to conclude the assertion of our main Theorem 2.2.

4.3. *Basic properties of the transport operator.* The following proposition holds.

PROPOSITION 4.4. *For any nonnegative* $F \in L^\infty(\mathcal{T}_D)$,

(4.6) $$\int \mathfrak{Q} F \, d\mathbb{P}_D = \int F \, d\mathbb{P}_D.$$

*Hence,* $\mathfrak{Q}$ *can be extended to a density-preserving operator* $\mathfrak{Q} : L^1(\mathcal{T}_D) \to L^1(\mathcal{T}_D)$.

PROOF. The left-hand side of (4.6) equals

$$\frac{1}{P_0[D = +\infty]} \iint \mathcal{K}(\omega, \omega') F(\omega) Q_{\omega'}[D = +\infty] \mathbb{P}(d\omega) \mathbb{P}(d\omega')$$

(4.7)
$$= \frac{1}{P_0[D = +\infty]}$$
$$\times \sum_{k,L=1}^{+\infty} \iiint \bar{\nu}_{S_k}(\pi; \omega, \omega') \mathbf{1}_{A(S_k)} \mathbf{1}_{[L-1,L)}(S_k)$$
$$\times Q_{T_{\pi(S_k)}\omega'}[D = +\infty] F(\omega) \mathbb{W}(d\pi) \mathbb{P}(d\omega) \mathbb{P}(d\omega').$$

Using properties (Z2) and (Z3) of operator $\mathcal{Z}$, we conclude that the right-hand side of (4.7) equals

$$\frac{1}{P_0[D=+\infty]} \sum_{k,L=1}^{+\infty} \int\int \nu_{S_k}(\pi;\omega)\mathbf{1}_{A(S_k)}\mathbf{1}_{[L-1,L)}(S_k)$$

$$\times Q_{T_{\pi(S_k)}\omega}[D=+\infty]F(\omega)\mathbb{W}(d\pi)\mathbb{P}(d\omega)$$

$$= \frac{1}{P_0[D=+\infty]}$$

$$\times \sum_{k=1}^{+\infty} \int \mathbf{M}_\omega[Q_{T_{\pi(S_k)}\omega}[D=+\infty], A(S_k), S_k < +\infty]F(\omega)\mathbb{P}(d\omega)$$

$$= \frac{1}{P_0[D=+\infty]} \int Q_\omega[D=+\infty, \tau_1 < +\infty]F(\omega)\mathbb{P}(d\omega).$$

Since $\tau_1 < +\infty$, $Q_u$-a.s. we conclude that the last expression is equal to the right-hand side of (4.6). $\square$

Suppose that $F:(\mathbb{R} \times \mathbb{R} \times \mathbb{R}^d)^\mathbb{N} \to \mathbb{R}$ and $G:\Omega \to \mathbb{R}$ are bounded and, respectively, Borel and $\mathcal{V}_0$-measurable. Let $n$, $N$ be positive integers, $0 \leq t_1 \leq \cdots \leq t_n$, and $F_1, \ldots, F_n:\mathbb{R}^d \to \mathbb{R}$, $H:(\mathbb{R} \times \mathbb{R} \times \mathbb{R}^d)^N \to \mathbb{R}$ be arbitrary bounded and measurable functions.

Define

$$(4.8) \qquad \xi_k := \int_{\tau_k}^{\tau_{k+1}} \prod_{p=1}^n F_p(\mathbf{u}(\pi(t_p+s)))\,ds$$

and

$$(4.9) \qquad \tilde{\xi}_k := (\xi_k, \tau_{k+1} - \tau_k, \pi(\tau_{k+1}) - \pi(\tau_k)).$$

Here $\tau_0 := 0$. Let also $q$ be a positive integer and $\pi^{(q)}(s) = \pi(s \wedge \tau_q)$,

$$(4.10) \qquad \xi_k^{(q)} := \int_{\tau_k}^{\tau_{k+1}} \prod_{p=1}^n F_p(\mathbf{u}(\pi^{(q)}(t_p+s)))\,ds.$$

and define $\tilde{\xi}_k^{(q)}$ accordingly.

PROPOSITION 4.5. *Let $n \geq 1$ be an arbitrary integer. Suppose that $0 \leq t_1 \leq \cdots \leq t_n$ are arbitrary and $\xi_k$, $\tilde{\xi}_k$, $k \geq 1$, are defined by formulas* (4.8) *and* (4.9). *Under the assumptions on $F$, $G$ specified above, we have:*

(i)

$$\iint F((\tilde{\xi}_{k+1})_{k\geq 1})G(\omega)P_D(d\omega,d\pi) \tag{4.11}$$
$$= \iint F((\tilde{\xi}_k)_{k\geq 1})\mathfrak{Q}G(\omega)P_D(d\omega,d\pi).$$

(ii) *In addition, suppose that $q \geq q_0 \geq N$ are certain integers, function $H$ and random variables $\xi_k^{(q)}$, $\tilde{\xi}_k^{(q)}$, $k \geq 1$, are specified as in the foregoing. Then, there exists a random variable $Y \in L^\infty(\mathcal{T}_D)$ such that*

$$\iint F((\tilde{\xi}_{k+q})_{k\geq 1})H(\tilde{\xi}_1^{(q_0)},\ldots,\tilde{\xi}_N^{(q_0)})G(\omega)P_D(d\omega,d\pi) \tag{4.12}$$
$$= \iint F((\tilde{\xi}_k)_{k\geq 1})\mathfrak{Q}^{q-q_0}Y(\omega)P_D(d\omega,d\pi).$$

$Y$ *is nonnegative when $G$, $H$ are nonnegative and*

$$(4.13) \quad \iint Y(\omega)P_D(d\omega,d\pi) = \iint H(\tilde{\xi}_1^{(q_0)},\ldots,\tilde{\xi}_N^{(q_0)})G(\omega)P_D(d\omega,d\pi).$$

PROOF. For any sequence $\mathbf{m} := (m_1,\ldots,m_q) \in \mathbb{Z}_+^q$, we define a sequence of Markovian times

$$(4.14) \quad \sigma_0^{\mathbf{m}} := 0 \quad \text{and} \quad \sigma_{r+1}^{\mathbf{m}} := \sigma_r^{\mathbf{m}} + S_{m_{r+1}} \circ \theta_{\sigma_r^{\mathbf{m}}}, \qquad r = 0,\ldots,q-1.$$

The sequence is defined on the set of paths satisfying

$$B(\mathbf{m}) := \Big[\pi \colon \text{all random times appearing in (4.14) are finite and}$$

$$\inf_{t\in[\sigma_r^{\mathbf{m}},\sigma_{r+1}^{\mathbf{m}}]}(\pi(t)\cdot\hat{\mathbf{v}} - \pi(\sigma_r^{\mathbf{m}})\cdot\hat{\mathbf{v}}) \geq -1, \forall\, r = 0,\ldots,q-1\Big].$$

Let

$$\tilde{\xi}_r^{\mathbf{m}} := \left(\int_{\sigma_r^{\mathbf{m}}}^{\sigma_{r+1}^{\mathbf{m}}} \prod_{p=1}^n F_p(\mathbf{u}(\pi(t_p+s)))\,ds, \sigma_{r+1}^{\mathbf{m}} - \sigma_r^{\mathbf{m}}, \pi(\sigma_{r+1}^{\mathbf{m}}) - \pi(\sigma_r^{\mathbf{m}})\right),$$
$$r = 0,\ldots,q-1,$$

and

$$\tilde{\xi}_r^{\mathbf{m},q_0} := \left(\int_{\sigma_r^{\mathbf{m}}}^{\sigma_{r+1}^{\mathbf{m}}} \prod_{p=1}^n F_p(\mathbf{u}(\pi^{(q_0)}(t_p+s)))\,ds, \sigma_{r+1}^{\mathbf{m}} - \sigma_r^{\mathbf{m}}, \pi(\sigma_{r+1}^{\mathbf{m}}) - \pi(\sigma_r^{\mathbf{m}})\right),$$
$$r = 0,\ldots,q-1.$$

We have

$$\frac{1}{P_0[D=+\infty]}\int\!\!\int F((\tilde{\xi}_{k+q})_{k\geq 1})H(\tilde{\xi}_1^{(q_0)},\ldots,\tilde{\xi}_N^{(q_0)})$$

$$\times G(\omega)\mathbf{1}_{[D(\pi)=+\infty]}P_0(d\omega,d\pi)$$

(4.15)
$$=\frac{1}{P_0[D=+\infty]}$$

$$\times \sum_{L=1}^{\infty}\sum_{\mathbf{m}}\int \mathbf{M}_\omega\Big[F((\tilde{\xi}_k\circ\theta_{\sigma_q^{\mathbf{m}}})_{k\geq 1})H(\tilde{\xi}_1^{\mathbf{m},q_0},\ldots,\tilde{\xi}_N^{\mathbf{m},q_0}),$$

$$D\circ\theta_{\sigma_q^{\mathbf{m}}}=+\infty, B(\mathbf{m}), L-1\leq\sigma_q^{\mathbf{m}}<L\Big]$$

$$\times G(\omega)\mathbb{P}(d\omega).$$

Using the strong Markov property and stationarity of the environment, we can recast the right-hand side of (4.15) in the form

(4.16) $$\frac{1}{P_0[D=+\infty]}\int\!\!\int F((\tilde{\xi}_{k+q-q_0})_{k\geq 1})Y(\omega)\mathbf{1}_{[D(\pi)=+\infty]}P_0(d\omega,d\pi),$$

where $Y(\omega)$ is a certain $\mathcal{V}_0$-measurable random variable. Note that $Y$ can be chosen so that it is nonnegative when $H$ and $G$ are nonnegative. Choosing $F\equiv 1$ in the argument above, we conclude also that $Y$ satisfies (4.13).

In the special case when $q=1$, $q_0=0$ and $H\equiv 1$, we can rewrite the right-hand side of (4.15) in the form

$$\frac{1}{P_0[D=+\infty]}\sum_{m,L=1}^{\infty}\int \mathbf{M}_\omega\Big[F((\tilde{\xi}_k\circ\theta_{S_m})_{k\geq 1}),$$

$$D\circ\theta_{S_m}=+\infty, A(S_m), L-1\leq S_m<L\Big]$$

$$\times G(\omega)\mathbb{P}(d\omega)$$

(4.17)
$$=\frac{1}{P_0[D=+\infty]}\sum_{m,L=1}^{\infty}\int \mathbf{M}_\omega\Big[\mathbf{M}_{\pi(S_m),\omega}[F((\tilde{\xi}_k)_{k\geq 1}), D=+\infty],$$

$$A(S_m), L-1\leq S_m<L\Big]$$

$$\times G(\omega)\mathbb{P}(d\omega).$$

Using Girsanov's theorem and subsequently conditioning on $[\pi(S_m) = \mathbf{x}, S_k \in [L-1, L)]$, we deduce that the right-hand side of (4.17) equals

$$\frac{1}{P_0[D = +\infty]} \sum_{m,L=1}^{\infty} \iint_{\mathbb{R}^d} \mathbf{M}_{m,L,\mathbf{x}}[\nu_{S_m}(\pi;\omega), A(S_m), L-1 \leq S_m < L]$$

(4.18)
$$\times \mathbf{M}_{T_\mathbf{x}\omega}[F((\tilde{\xi}_k)_{k\geq 1}), D = +\infty]$$

$$\times G(\omega)\mathbb{P}(d\omega)\mathbb{F}_m(d\mathbf{x}).$$

Since the second and third factors of the integrand appearing in (4.18) are, respectively, $\mathcal{R}$ and $\mathcal{V}_0$-measurable, we can rewrite the entire expression, using property (Z3) of operator $\mathcal{Z}$, in the following form [cf. (4.5)]

$$\frac{1}{P_0[D = +\infty]} \sum_{m,L=1}^{\infty} \iiint_{\mathbb{R}^d} \mathcal{M}_{m,L,\mathbf{x}}(\omega, \omega')$$

$$\times \mathbf{M}_{T_\mathbf{x}\omega'}[F((\tilde{\xi}_k)_{k\geq 1}), D = +\infty]$$

(4.19)
$$\times G(\omega)\mathbb{P}(d\omega)\mathbb{P}(d\omega')\mathbb{F}_m(d\mathbf{x})$$

$$= \int \mathbf{M}_\omega[F((\tilde{\xi}_k)_{k\geq 1}), D = +\infty]\mathfrak{Q}G(\omega)\mathbb{P}(d\omega).$$

We have proved therefore (4.11). To obtain (4.12), thus finishing the proof of the proposition, it suffices only to apply (4.17) to (4.19) $q - q_0$ times. □

Suppose that $m, n \geq 1$, $0 \leq t_1 \leq \cdots \leq t_n$, $0 \leq s_1 \leq \cdots \leq s_m$ and $F: (\mathbb{R}^d)^n \to \mathbb{R}$, $G: (\mathbb{R}^d)^m \to \mathbb{R}$. Let

$$\mathcal{U}_k := F(\mathbf{u}(\pi(t_1 + \tau_k)), \ldots, \mathbf{u}(\pi(t_n + \tau_k)))$$

and

$$\mathcal{G}^{(q_0)} := G(\mathbf{u}(\pi^{(q_0)}(s_1)), \ldots, \mathbf{u}(\pi^{(q_0)}(s_m))).$$

The proof of the proposition formulated below is analogous to the one used to show Proposition 4.5.

PROPOSITION 4.6. *Suppose that we are given $q \geq q_0$. Then, there exists a random variable $Y$ such that $Y \in L^\infty(\mathcal{T}_D)$ and*

$$\iint \mathcal{U}_q(\omega, \pi)\mathcal{G}^{(q_0)}(\omega, \pi)P_D(d\omega, d\pi)$$

(4.20)
$$= \iint \mathcal{U}_0(\omega, \pi) \mathfrak{Q}^{q-q_0} Y(\omega) P_D(d\omega, d\pi).$$

$Y$ is nonnegative when $\mathcal{G}^{(q_0)}$ is nonnegative and

(4.21) $$\iint Y(\omega) P_D(d\omega, d\pi) = \iint \mathcal{G}^{(q_0)}(\omega, \pi) P_D(d\omega, d\pi).$$

4.4. *The existence of an invariant density.* The following result is of crucial importance for us in the sequel.

THEOREM 4.7. *There exists a unique $H_* \in \mathcal{D}(\mathcal{T}_D)$ such that $\mathfrak{Q} H_* = H_*$ and $H_* > 0$, $\mathbb{P}_D$-a.s. (thus also $\mathbb{P}$-a.s.). In addition, there exist deterministic constants $\gamma_4 \in (0,1)$, $\gamma_5 > 0$, such that*

(4.22) $$\int |\mathfrak{Q}^n F - H_*| \, d\mathbb{P}_D \leq \gamma_5 \gamma_4^n \qquad \forall F \in \mathcal{D}(\mathcal{T}_D), n \geq 1.$$

The existence and uniqueness of a positive invariant density is a consequence of Theorem 5.6.2 of [11] and the following lemma.

LEMMA 4.8. *There exists a deterministic constant $\gamma_3 > 0$ such that $\mathfrak{Q} F \geq \gamma_3$, $\mathbb{P}_D$-a.s. for all $F \in \mathcal{D}(\mathcal{T}_D)$.*

PROOF. Suppose that $A \in \mathcal{V}_0$. We have

$$\int_A \mathfrak{Q} F \, d\mathbb{P}_D = \frac{1}{P_0[D = \infty]}$$

(4.23)
$$\times \sum_{k=1}^{+\infty} \iint \mathbf{M}_{\omega,\omega'}[Q_{\pi(S_k),\omega'}[D = +\infty]$$

$$\times \mathbf{1}_A(T_{\pi(S_k)}\omega'), A(S_k), S_k < +\infty]$$

$$\times F(\omega) \mathbb{P}(d\omega) \mathbb{P}(d\omega').$$

Here $\mathbf{M}_{\omega,\omega'}$ is the expectation operator corresponding to $Q_{\omega,\omega'}$—the unique solution of the martingale problem

(4.24)
$$d\mathbf{x}(t) = \mathbf{U}(\mathbf{x}(t); \omega, \omega') \, dt + d\mathbf{w}(t),$$
$$\mathbf{x}(0) = \mathbf{0}.$$

Using Lemma 3.1, we can estimate $Q_{\pi(S_1),\omega'}[D = +\infty] \geq \gamma$, while $A(S_1)$ [no backtracking may occur before $S_1$; cf. (3.4)] is clearly implied by the event $[S_1 \leq D]$. Hence the right-hand side of (4.23) is bound from below by

$$(4.25) \quad \gamma \int\int \mathbf{M}_{\omega,\omega'}[\mathbf{1}_A(T_{\pi(S_1)}\omega'), S_1 \leq D, S_1 < +\infty] F(\omega) \mathbb{P}(d\omega) \mathbb{P}(d\omega').$$

Let $G$ be a certain bounded subregion of the layer $[\mathbf{x} \in \mathbb{R}^d : -1 \leq \mathbf{x} \cdot \hat{\mathbf{v}} \leq r_0 + 1]$ containing $\mathbf{0}$, with a $C^\infty$-smooth boundary $\partial G$. We assume further that part of $\partial G$ of positive surface measure $S$ is contained in the hyperplane $\mathbb{H} := [\mathbf{x} \in \mathbb{R}^d : \mathbf{x} \cdot \hat{\mathbf{v}} = r_0 + 1]$. The expression in (4.25) can be further estimated from below by

$$(4.26) \quad \gamma \int\int \mathbf{M}_{\omega,\omega'}[\varphi(\pi(\tau_G); \omega')] F(\omega) \mathbb{P}(d\omega) \mathbb{P}(d\omega').$$

Here $\tau_G$ denotes the exit time from $G$, $\varphi(\mathbf{x}; \omega') := \mathbf{1}_A(T_\mathbf{x}\omega')$ for $\mathbf{x} \in \partial G \cap \mathbb{H}$ and $\varphi(\mathbf{x}; \omega') := 0$, if $\mathbf{x} \in \partial G \setminus \mathbb{H}$. Using absolute continuity of the harmonic measure and the standard lower bounds for Green's function corresponding to the generator of (4.24) and the region $G$ (see, e.g., [4], Theorem 3.1, page 616), we can bound (4.26) from below by

$$(4.27) \quad \begin{aligned} & C \int\int \left[ \int_{\partial G \cap \mathbb{H}} \mathbf{1}_A(T_\mathbf{y}\omega') S(d\mathbf{y}) \right] F(\omega) \mathbb{P}(d\omega) \mathbb{P}(d\omega') \\ &= C |\partial G \cap \mathbb{H}| \mathbb{P}[A] \int F(\omega) \mathbb{P}(d\omega) \\ &\geq C P_0[D = +\infty] |\partial G \cap \mathbb{H}| \mathbb{P}[A] \int F(\omega) \mathbb{P}_D(d\omega) \\ &= C P_0[D = +\infty] |\partial G \cap \mathbb{H}| \mathbb{P}[A], \end{aligned}$$

where $C > 0$ is a certain deterministic, positive constant and $|E|$ denotes the surface measure of a Lebesgue measurable subset $E \subset \mathbb{H}$. Here we used the fact that $F \in \mathcal{D}(\mathcal{T}_D)$. $\square$

**5. The construction of an invariant measure.** Let $H_*$ be the invariant density for $\mathfrak{Q}$, see Theorem 4.7. We set

$$P_{H_*}(d\omega, d\pi) := H_*(\omega) P_D(d\omega, d\pi).$$

Note that $P_{H_*}$ is a probability measure on $(\Omega \times \mathfrak{X}, \overline{\mathcal{B}(\Omega)} \otimes \mathcal{M})$.

THEOREM 5.1. *Let $n \geq 1$ be an integer and $0 \leq t_1 \leq \cdots \leq t_n$ be arbitrary. We suppose that $F_1, \ldots, F_n : \mathbb{R}^d \to \mathbb{R}$ are any bounded measurable functions and the sequence $(\tilde{\xi}_k)_{k \geq 0}$ is given by (4.8) and (4.9). Then $(\tilde{\xi}_k)_{k \geq 0}$ is stationary and ergodic over the probability space $(\Omega \times \mathfrak{X}, \overline{\mathcal{B}(\Omega)} \otimes \mathcal{M}, P_{H_*})$.*

PROOF. Stationarity is a direct consequence of part (i) of Proposition 4.5 and the definition of $H_*$. To prove ergodicity, we show that any bounded measurable function $F:(\mathbb{R} \times \mathbb{R} \times \mathbb{R}^d)^{\mathbb{N}} \to \mathbb{R}$, for which

(5.1) $$F((\tilde{\xi}_{k+n})_{k\geq 1}) = F((\tilde{\xi}_k)_{k\geq 1}) \qquad \forall n \geq 1, \ P_{H_*}\text{-a.s.},$$

satisfies $F((\tilde{\xi}_k)_{k\geq 1}) \equiv \text{const}$, $P_{H_*}$-a.s. Let $\varepsilon > 0, N \geq 1$, be arbitrary. We can find $F^{(N)}:(\mathbb{R} \times \mathbb{R} \times \mathbb{R}^d)^N \to \mathbb{R}$ bounded, continuous and such that

$$\iint |F((\tilde{\xi}_k)_{k\geq 1}) - F^{(N)}(\tilde{\xi}_1, \ldots, \tilde{\xi}_N)| dP_{H_*} < \varepsilon.$$

Then,

(5.2) $$\iint |F((\tilde{\xi}_k)_{k\geq 1})[F((\tilde{\xi}_k)_{k\geq 1}) - F^{(N)}(\tilde{\xi}_1, \ldots, \tilde{\xi}_N)]| \, dP_{H_*} < \varepsilon \sup |F|.$$

On the other hand, for any $q \geq q_0$, we have, from (5.1),

(5.3) $$\iint F((\tilde{\xi}_k)_{k\geq 1}) F^{(N)}(\tilde{\xi}_1^{(q_0)}, \ldots, \tilde{\xi}_N^{(q_0)}) \, dP_{H_*}$$
$$= \iint F((\tilde{\xi}_{k+q})_{k\geq 1}) F^{(N)}(\tilde{\xi}_1^{(q_0)}, \ldots, \tilde{\xi}_N^{(q_0)}) \, dP_{H_*}.$$

By virtue of Proposition 4.5, we conclude that the right-hand side of (5.3) equals

$$\iint F((\tilde{\xi}_k)_{k\geq 1}) \mathfrak{Q}^{q-q_0} Y \, dP_{H_*}$$

for a certain $\mathcal{V}_0$-measurable $Y$ such that

$$\iint Y \, dP_{H_*} = \iint F^{(N)}(\tilde{\xi}_1^{(q_0)}, \ldots, \tilde{\xi}_N^{(q_0)}) \, dP_{H_*}.$$

First letting $q \uparrow +\infty$ and then $q_0 \uparrow +\infty$, we conclude that

(5.4) $$\iint F((\tilde{\xi}_k)_{k\geq 1}) F^{(N)}(\tilde{\xi}_1, \ldots, \tilde{\xi}_N) \, dP_{H_*}$$
$$= \iint F((\tilde{\xi}_k)_{k\geq 1}) \, dP_{H_*} \iint F^{(N)}(\tilde{\xi}_1, \ldots, \tilde{\xi}_N) \, dP_{H_*},$$

which, in light of (5.2), yields

$$\left| \iint [F((\tilde{\xi}_k)_{k\geq 1})]^2 \, dP_{H_*} - \left[ \iint F((\tilde{\xi}_k)_{k\geq 1}) \, dP_{H_*} \right]^2 \right| < 2\varepsilon \sup |F|.$$

Since $\varepsilon > 0$ was chosen arbitrarily, we conclude that $F((\tilde{\xi}_k)_{k\geq 1}) \equiv \text{const}$, $P_{H_*}$-a.s. $\square$

The following proposition holds.

PROPOSITION 5.2. *We have*

(5.5) $$\iint \tau_1 \, dP_{H_*} < +\infty$$

*and*

(5.6) $$\iint |\pi(\tau_1)| \, dP_{H_*} < +\infty.$$

PROOF. First we show the following.

LEMMA 5.3.

(5.7) $$w_* := \iint \hat{\mathbf{v}} \cdot \pi(\tau_1) \, dP_{H_*} < +\infty.$$

PROOF. We can write that

$$\hat{\mathbf{v}} \cdot \pi(\tau_1) = \sum_{k=0}^{K-1} \hat{\mathbf{v}} \cdot (\pi(S_{k+1}) - \pi(S_k))$$

with random variable $K$ defined before the statement of Corollary 3.3. Hence,

(5.8) $$\hat{\mathbf{v}} \cdot \pi(\tau_1) \leq r_0 + 1 + \sum_{k=1}^{K-1} (r_0 + 1 + M_k - \hat{\mathbf{v}} \cdot \pi(S_k)),$$

and in consequence,

$$\iint \hat{\mathbf{v}} \cdot \pi(\tau_1) \, dP_{H_*} \leq r_0 + 1 + \sum_{1 \leq k' < k} \iint (r_0 + 1 + M_{k'} - \hat{\mathbf{v}} \cdot \pi(S_{k'}))$$

(5.9) $$\times \mathbf{1}_{[S_k < +\infty, D \circ \theta_{S_k} = +\infty]} \, dP_{H_*}$$

$$\leq r_0 + 1 + \sum_{1 \leq k' < k} \iint (r_0 + 1 + M_{k'} - \hat{\mathbf{v}} \cdot \pi(S_{k'}))$$

$$\times \mathbf{1}_{[R_{k-1} < +\infty, D \circ \theta_{S_k} = +\infty]} \, dP_{H_*}.$$

Since $R_{k-1} = D \circ \theta_{S_{k-1}} + S_{k-1}$, we obtain, upon a multiple application of the strong Markov property for $Q_\omega$ and (3.7), that the right-hand side of (5.9) is less than or equal to

$$r_0 + 1 + \sum_{1 \leq k' < k} (1 - \gamma)^{k-1-k'}$$

$$\times \iint (r_0 + 1 + M_{k'} - \hat{\mathbf{v}} \cdot \pi(S_{k'})) \mathbf{1}_{[R_{k'} < +\infty]} \, dP_{H_*}$$

$$(5.10) \quad \leq r_0 + 1 + \frac{1}{P_{\mathbf{0}}[D = +\infty]} \sum_{1 \leq k' < k} (1 - \gamma)^{k-1-k'}$$

$$\times \int H_*(\omega) \mathbf{M}_\omega [\mathbf{M}_{\pi(S_{k'}),\omega}[r_0 + 1 + M_*, D < +\infty],$$

$$S_{k'} < +\infty] \mathbb{P}(d\omega).$$

By virtue of (3.9), we conclude that the right-hand side of (5.10) is less than or equal to

$$r_0 + 1 + \frac{r_0 + 1 + \gamma_2}{P_{\mathbf{0}}[D = +\infty]}$$

$$(5.11) \quad \times \sum_{1 \leq k' < k} (1 - \gamma)^{k-1-k'} \int H_*(\omega) Q_\omega[S_{k'} < +\infty] \mathbb{P}(d\omega)$$

$$\stackrel{(3.11)}{\leq} r_0 + 1 + C \sum_{k=1}^{+\infty} k(1 - \gamma)^{k-1} < +\infty,$$

for some constant $C > 0$, and (5.7) follows. We have used here the fact that

$$\int H_*(\omega) \mathbb{P}(d\omega) \stackrel{\text{Lemma 3.1}}{\leq} \frac{P_{\mathbf{0}}[D = +\infty]}{\gamma} \int H_*(\omega) \mathbb{P}_D(d\omega) < +\infty. \quad \square$$

Continuing with the proof of the proposition, we let $(k_m)_{m \geq 1}$ be a random sequence of integers defined by

$$(5.12) \quad \tau_{k_m} \leq U_m < \tau_{k_m+1}.$$

Recall the convention that $\tau_0 := 0$. Then $P_{H_*}$-a.s. we have

$$\hat{\mathbf{v}} \cdot \pi(\tau_{k_m}) \leq m < \hat{\mathbf{v}} \cdot \pi(\tau_{k_m+1}).$$

By virtue of Theorem 5.1 and the individual ergodic theorem, we conclude that

$$\lim_{k \uparrow +\infty} \frac{\hat{\mathbf{v}} \cdot \pi(\tau_k)}{k} = w_* \stackrel{\text{Lemma 5.3}}{<} +\infty, \qquad P_{H_*}\text{-a.s.}$$

But

$$\frac{\hat{\mathbf{v}} \cdot \pi(\tau_{k_m})}{k_m} \leq \frac{m}{k_m} < \frac{\hat{\mathbf{v}} \cdot \pi(\tau_{k_m+1})}{k_m},$$

therefore

$$\lim_{m \uparrow +\infty} \frac{k_m}{m} = \frac{1}{w_*}, \qquad P_{H_*}\text{-a.s.}$$

Let

$$t_* := \iint \tau_1 \, dP_{H_*}. \tag{5.13}$$

Trivially, we conclude that $t_* \in (0, +\infty]$. We also have

$$\frac{U_m}{m} \geq \frac{\tau_{k_m}}{k_m} \times \frac{k_m}{m}. \tag{5.14}$$

Inequality (3.10) of Lemma 3.2 implies in particular that

$$\liminf_{m \to +\infty} \frac{U_m}{m} < +\infty, \qquad P_{H_*}\text{-a.s.} \tag{5.15}$$

On the other hand, an application of the ergodic theorem to the sequence $((\tau_{n+1} - \tau_n))_{n \geq 1}$ implies that the right-hand side of (5.14) tends $P_{H_*}$-a.s. to $t_*/w_*$, which by virtue of Lemma 5.3 and (5.15) belongs to $(0, +\infty)$. Consequently, we conclude that $t_* < +\infty$ and (5.5) holds.

Additionally, we have

$$\iint |\pi(\tau_1)| \, dP_{H_*} \leq \|\mathbf{u}\|_{L^\infty} \iint \tau_1 \, dP_{H_*} + \iint |\mathbf{w}(\tau_1)| \, dP_{H_*}. \tag{5.16}$$

Denoting $X := \sup_{0 \leq t \leq 1} |\mathbf{w}(t)|$, $Y := \sup_{t \geq 1} |\mathbf{w}(t)| t^{-3/4}$, we can estimate the second term on the right-hand side of (5.16) by

$$\iint X \mathbf{1}_{[\tau_1 \leq 1]} \, dP_{H_*} + \iint Y \tau_1^{3/4} \, dP_{H_*}$$
$$\leq \iint X \, dP_{H_*} + \left(\iint Y^4 \, dP_{H_*}\right)^{1/4} \left(\iint \tau_1 \, dP_{H_*}\right)^{3/4} < +\infty.$$

Here we used the fact that $Q_\omega[Y > u] \leq c_1 \exp\{-c_2 u^2\}$ for some constants $c_1, c_2 > 0$ independent of $\omega$ and all $u > 0$. This can easily be concluded from, for example, [1], Theorem 5.2, page 120. □

As a consequence of Proposition 5.2, Theorem 5.1 and the individual ergodic theorem, we obtain the following.

COROLLARY 5.4. *We have*

$$\frac{1}{N} \sum_{k=0}^{N-1} \xi_k \to \iint \left[ \int_0^{\tau_1} \prod_{p=1}^n F_p(\mathbf{u}(\pi(t_p + s))) \, ds \right] P_{H_*}(d\omega, d\pi) \tag{5.17}$$

$$\text{as } N \to +\infty.$$

*The convergence in* (5.17) *holds both* $P_{H_*}$*-a.s. and in the* $L^1(P_{H_*})$*-sense.*

Suppose that $Q_{\mathbf{x},\omega,t}$ is the probability measure on $\mathfrak{X}_t$ obtained by restricting $Q_{\mathbf{x},\omega}$ to $\mathcal{M}_t$. Let also $Q_{\omega,t}^{\mathbf{x},\mathbf{y}}$, $\mathbf{M}_{\omega,t}^{\mathbf{x},\mathbf{y}}$ denote, respectively, the probability measure and the respective expectation obtained by conditioning $Q_{\mathbf{x},\omega,t}$ on the event $\pi(t) = \mathbf{y}$. Recall (from Section 2) that $p^\omega(\cdot,\cdot;\cdot,\cdot)$ denotes the transition probability density of the diffusion given by (1.1). We set

$$\mathcal{H}_m(\mathbf{x}, s, \omega, \pi)$$
$$:= \mathbf{1}_{[D(\mathbf{x}\cdot\hat{\mathbf{v}})=+\infty]}(\pi) p^\omega(s, \mathbf{x}, \mathbf{0}) Q_{\omega,s}^{\mathbf{x},\mathbf{0}}[A(s), S_m \leq s < S_{m+1}] H_*(T_\mathbf{x}\omega).$$

LEMMA 5.5. *Let* $n \geq 1$, $F_1, \ldots, F_n \in C_b(\mathbb{R}^d)$ *and* $0 \leq t_1 \leq \cdots \leq t_n$. *Then,*

(5.18)
$$\sum_{m=1}^{+\infty} \int_0^{+\infty} \int_{\mathbb{R}^d} \int\int \prod_{p=1}^n F_p(\mathbf{u}(\pi(t_p))) \mathcal{H}_m(\mathbf{x}, s, \omega, \pi) \, ds \, d\mathbf{x} \, P_0(d\omega, d\pi)$$
$$= \int \mathbf{M}_\omega\left[\int_0^{\tau_1} \prod_{p=1}^n F_p(\mathbf{u}(\pi(t_p + s))) \, ds, D = +\infty\right] H_*(\omega) \mathbb{P}(d\omega).$$

REMARK 5.6. Note that in light of Proposition 5.2, the right-hand side of (5.18) is finite.

PROOF OF LEMMA 5.5. Before proceeding with the proof, we introduce two additional renewal structures via a slight modification of the times $(S_k)_{k \geq 1}$. These structures allow us to describe moments $S_k$ that occur after certain deterministically fixed time $s$; see (5.26).

Let $l, m \in \mathbb{R}$. Recall that $D(l)$ is defined as in (3.1). Let $M_0(l) := \max[\pi(t) \cdot \hat{\mathbf{v}} : 0 \leq t \leq D(l)]$. We define the stopping time $S_1^{(1)}(l, m)$ as follows. On the event $D(l) < +\infty$, we let

$$S_1^{(1)}(l, m) := \min[t \geq D(l) : \pi(t) \cdot \hat{\mathbf{v}} \geq (M_0(l) \vee m) + r_0 + 1]$$

and

$$R_1^{(1)}(l, m) := D \circ \theta_{S_1^{(1)}(l,m)} + S_1^{(1)}(l, m),$$

$$M_1^{(1)}(l, m) := \max[\pi(t) \cdot \hat{\mathbf{v}} : 0 \leq t \leq R_1^{(1)}(l, m)].$$

We set

$$S_1^{(1)}(l, m) := R_1^{(1)}(l, m) := +\infty \qquad \text{if } D(l) = +\infty.$$

The subsequent times $R_k^{(1)}(l, m), S_k^{(1)}(l, m)$ and maxima $M_k^{(1)}(l, m)$ are defined as follows:

(5.19) $\qquad S_{k+1}^{(1)}(l, m) := U_{M_k(l,m)+r_0+1},$

(5.20) $\qquad R_{k+1}^{(1)}(l, m) := S_{k+1}^{(1)}(l, m) + D \circ \theta_{S_{k+1(l,m)}^{(1)}},$

(5.21) $\qquad M_{k+1}^{(1)}(l, m) := \max[\pi(t) \cdot \hat{\mathbf{v}} : 0 \leq t \leq R_k^{(1)}(l, m)].$

Similarly for $l \geq \pi(0) \cdot \hat{\mathbf{v}}$, we define

$$S_1^{(2)}(l) := \min[t : \pi(t) \cdot \hat{\mathbf{v}} \geq l + r_0 + 1]$$

and

$$R_1^{(2)}(l) := D \circ \theta_{S_1^{(2)}(l)} + S_1^{(2)}(l),$$

$$M_1^{(2)}(l) := \max[\pi(t) \cdot \hat{\mathbf{v}} : 0 \leq t \leq R_1^{(2)}(l)].$$

The subsequent times $R_k^{(2)}(l), S_k^{(2)}(l)$ and maxima $M_k^{(2)}(l)$ are defined by means of (5.19)–(5.21) with the obvious replacement of superscripts and arguments $(l, m)$ by $l$. Let

$$K^{(1)}(l, m) := \min[k : S_k^{(1)}(l, m) < +\infty, D \circ \theta_{S_k^{(1)}(l,m)} = +\infty]$$

and

$$K^{(2)}(l) := \min[k : S_k^{(2)}(l) < +\infty, D \circ \theta_{S_k^{(2)}(l)} = +\infty].$$

Note that the definitions of $S_k^{(1)}(\cdot, \cdot)$, $k \geq 1$, differ from the respective definitions of $S_k$, $k \geq 1$, only through the designation of the first time $S_1^{(1)}(\cdot, \cdot)$. The same remark extends also to $S_k^{(2)}(\cdot)$, $k \geq 1$. Therefore, a straightforward adaptation of the argument used to prove Corollary 3.3 allows us to conclude that, for each $l, m \in \mathbb{R}$, $\mathbf{x} \in \mathbb{R}^d$, we have

(5.22) $\quad Q_{\mathbf{x},\omega}[K^{(1)}(l, m) < +\infty, S_{K^{(1)}(l,m)}^{(1)}(l, m) < +\infty] = 1, \qquad \mathbb{P}\text{-a.s.}$,

(5.23) $\qquad Q_{\mathbf{x},\omega}[K^{(2)}(l) < +\infty, S_{K^{(2)}(l)}^{(2)}(l) < +\infty] = 1, \qquad \mathbb{P}\text{-a.s.}$

To explain the meaning of $S_k^{(1)}(\cdot, \cdot)$, $S_k^{(2)}(\cdot)$, $k \geq 1$, consider the case when $s$ is a certain fixed deterministic time, $m_1 < m_2$ are two positive integers and $\pi(\cdot)$ is a path that satisfies $s \in [S_{m_1}, R_{m_1}), S_{m_2} < +\infty$. Then, we can write

(5.24) $\qquad S_{m_2}(\pi) = s + S_{m_2-m_1}^{(1)}(\pi(S_{m_1}) \cdot \hat{\mathbf{v}}, N_{m_1}(s)) \circ \theta_s(\pi).$

Here $N_m(s) := \max[\pi(t) \cdot \hat{\mathbf{v}} : t \in [S_m \wedge s, s]]$. If, on the other hand, $s \in [R_{m_1}, S_{m_1+1})$, $S_{m_2} < +\infty$, we have

(5.25) $\qquad S_{m_2}(\pi) = s + S_{m_2-m_1}^{(2)}(M_{m_1}(\pi)) \circ \theta_s(\pi).$

Recall that $M_{m_1}$ is defined in (3.5) and (3.6).

For brevity let us denote

$$F(s) := \prod_{p=1}^{n} F_p(\mathbf{u}(\theta_s(\pi)(t_p)))$$

and $B_m := [D \circ \theta_{S_m} = +\infty, D = +\infty]$. The right-hand side of (5.18) equals

$$\sum_{m=1}^{+\infty} \int \mathbf{M}_\omega \left[ \int_0^{S_m} F(s) \, ds, S_m < +\infty, B_m \right] H_*(\omega) \mathbb{P}(d\omega)$$

$$= \sum_{0 \leq m_1 \leq m_2 - 1} \int_0^{+\infty} \left\{ \int \mathbf{M}_\omega [F(s) \mathbf{1}_{[S_{m_1}, R_{m_1})}(s), \right.$$

(5.26)
$$\left. S_{m_2} < +\infty, B_{m_2}] H_*(\omega) \mathbb{P}(d\omega) \right\} ds$$

$$+ \sum_{0 \leq m_1 \leq m_2 - 1} \int_0^{+\infty} \left\{ \int \mathbf{M}_\omega [F(s) \mathbf{1}_{[R_{m_1}, S_{m_1+1})}(s), \right.$$

$$\left. S_{m_2} < +\infty, B_{m_2}] H_*(\omega) \mathbb{P}(d\omega) \right\} ds.$$

Using the Markov property, the definition of stopping times $S_k^{(1)}(l, m)$ and (5.24) we can recast the first term on the right-hand side of (5.26) as being equal to [cf. (3.4)]

$$\sum_{0 \leq m_1 \leq m_2 - 1} \int_0^{+\infty} \left\{ \int \mathbf{M}_\omega [\mathbf{1}_{[S_{m_1}, R_{m_1})}(s) \right.$$

(5.27)
$$\times g_{m_2 - m_1}(\pi(s), \pi(S_{m_1}) \cdot \hat{\mathbf{v}}, N_{m_1}(s)), A(s)]$$

$$\left. \times H_*(\omega) \mathbb{P}(d\omega) \right\} ds,$$

where

$$g_k(\mathbf{x}, l, m) := \mathbf{M}_{\mathbf{x}, \omega}[F(0), S_k^{(1)}(l, m) < +\infty, D \circ \theta_{S_k^{(1)}(l, m)} = +\infty, D(0) = +\infty].$$

Using (5.22), we conclude that the expression in (5.27) equals

$$\sum_{m_1 = 0}^{+\infty} \int_0^{+\infty} \left\{ \int \mathbf{M}_\omega [\mathbf{1}_{[S_{m_1}, R_{m_1})}(s) A(s)] H_*(\omega) \right.$$

(5.28)
$$\left. \times \mathbf{M}_{\pi(s), \omega}[F(0), D(0) = +\infty] \mathbb{P}(d\omega) \right\} ds.$$

Conditioning on the event $\pi(s) = \mathbf{x}$, we obtain that the expression in (5.28) equals

$$\sum_{m_1 = 0}^{+\infty} \int \int_0^{+\infty} \int_{\mathbb{R}^d} p^\omega(s, \mathbf{0}, \mathbf{x}) \mathbf{M}_{\omega, s}^{\mathbf{0}, \mathbf{x}}[\mathbf{1}_{[S_{m_1}, R_{m_1})}(s), A(s)]$$

(5.29)
$$\times \mathbf{M}_{\mathbf{x}, \omega}[F(0), D(0) = +\infty] H_*(\omega) \mathbb{P}(d\omega) \, ds \, d\mathbf{x}.$$

Using the homogeneity of $\mathbb{P}$ and changing variables $\mathbf{x} := -\mathbf{x}$, we conclude that the expression in (5.29) equals

$$\sum_{m_1=0}^{+\infty} \iint \int_0^{+\infty} \int_{\mathbb{R}^d} \mathbf{1}_{[D(\mathbf{x}\cdot\hat{\mathbf{v}})=+\infty]}(\pi) p^\omega(s,\mathbf{x},\mathbf{0})$$

(5.30)
$$\times \mathbf{M}_{\omega,s}^{\mathbf{x},\mathbf{0}}[A(s), S_{m_1} \leq s < R_{m_1}]$$

$$\times H_*(T_\mathbf{x}\omega) F(0) P_0(d\omega, d\pi) \, ds \, d\mathbf{x}.$$

Note that thanks to the definition of $A(s)$, the integration over variable $\mathbf{x}$ extends only over the region $[\mathbf{x} : \mathbf{x} \cdot \hat{\mathbf{v}} \leq 1]$.

Repeating the same type of calculations for the second term on the right-hand side of (5.26) [using stopping times $S_k^{(2)}(l)$ instead of $S_k^{(1)}(l,m)$ and (5.25)], we conclude that it equals

$$\sum_{m_1=0}^{+\infty} \iint \int_0^{+\infty} \int_{\mathbb{R}^d} \mathbf{1}_{[D(\mathbf{x}\cdot\hat{\mathbf{v}})=+\infty]}(\pi) p^\omega(s,\mathbf{x},\mathbf{0})$$

(5.31)
$$\times \mathbf{M}_{\omega,s}^{\mathbf{x},\mathbf{0}}[A(s), R_{m_1} \leq s < S_{m_1+1}]$$

$$\times H_*(T_\mathbf{x}\omega) F(0) P_0(d\omega, d\pi) \, ds \, d\mathbf{x},$$

and (5.18) follows. □

Applying Lemma 5.5 to $n = 1$, $t_1 = 1$ and $F_1 \equiv 1$, we conclude immediately the following.

COROLLARY 5.7.

(5.32)
$$P_0[D = +\infty] \iint \tau_1 P_{H_*}(d\omega, d\pi)$$
$$= \sum_{m=1}^{+\infty} \int_0^{+\infty} \int_{\mathbb{R}^d} \iint \mathcal{H}_m(\mathbf{x}, s, \omega, \pi) \, ds \, d\mathbf{x} \, P_0(d\omega, d\pi).$$

Let

(5.33) $$\mu(d\omega, d\pi) := \frac{1}{Z} \sum_{m=1}^{+\infty} \left[ \int_0^{+\infty} \int_{\mathbb{R}^d} \mathcal{H}_m(\mathbf{x}, s, \omega, \pi) \, ds \, d\mathbf{x} \right] P_0(d\omega, d\pi),$$

where the constant $Z$, by definition, equals the right-hand side of (5.32). Thanks to (5.5), we have $Z < +\infty$. By virtue of Corollary 5.7 $\mu$ is a probability measure.

PROPOSITION 5.8. *The process $V(\cdot)$ given over $(\Omega \times \mathfrak{X}, \overline{\mathcal{B}(\Omega)} \otimes \mathcal{M}, \mu)$ by formula (2.3) is stationary and ergodic.*

PROOF. *The proof of stationarity.* Let $n \geq 1$, $F_1, \ldots, F_n \in C_b(\mathbb{R}^d)$ and $0 \leq t_1 \leq \cdots \leq t_n$. Then, for any $h \geq 0$, we can write

$$\iint \prod_{p=1}^n F_p(\mathbf{u}(\pi(t_p + h))) \mu(d\omega, d\pi)$$

$$\stackrel{\text{Lemma 5.5}}{=} \frac{P_0[D = +\infty]}{Z} \iint \left[ \int_0^{\tau_1} \prod_{p=1}^n F_p(\mathbf{u}(\pi(t_p + h + s))) \, ds \right]$$

$$\times P_{H_*}(d\omega, d\pi)$$

(5.34)

$$\stackrel{\text{Corollary 5.4}}{=} \frac{P_0[D = +\infty]}{Z}$$

$$\times \lim_{N\uparrow+\infty} \frac{1}{N} \iint \left[ \int_0^{\tau_N} \prod_{p=1}^n F_p(\mathbf{u}(\pi(t_p + h + s))) \, ds \right]$$

$$\times P_{H_*}(d\omega, d\pi).$$

Since the integration over an interval of length $h$ does not influence the value of the expression on the utmost right-hand side of (5.34), we conclude that it is in fact equal to

$$\frac{P_0[D = +\infty]}{Z} \lim_{N\uparrow+\infty} \frac{1}{N} \iint \left[ \int_0^{\tau_N} \prod_{p=1}^n F_p(\mathbf{u}(\pi(t_p + s))) \, ds \right] P_{H_*}(d\omega, d\pi)$$

$$= \iint \prod_{p=1}^n F_p(\mathbf{u}(\pi(t_p))) \mu(d\omega, d\pi).$$

*Proof of ergodicity.* Proving ergodicity is tantamount to showing that, for any bounded and Borel measurable $F : \mathfrak{X} \to \mathbb{R}$ that satisfies

(5.35) $\qquad F \circ \theta_t(V(\cdot)) = F(V(\cdot)) \qquad \forall t \geq 0, \mu\text{-a.s.},$

we have $F(V(\cdot)) = \text{const}$, $\mu$-a.s. Similarly to what we have done in the proof of Theorem 5.1, for any $\varepsilon > 0$, we can find $N \geq 1$, $0 \leq t_1 \leq \cdots \leq t_N$ and $F^{(N)} : (\mathbb{R}^d)^N \to \mathbb{R}$ bounded, continuous such that

(5.36) $\qquad \iint \left| F(V(\cdot)) - F^{(N)}(V(t_1), \ldots, V(t_N)) \right| d\mu < \varepsilon.$

Let $q \geq q_0$ be arbitrary and $V^{(q_0)}(t) := V(t \wedge \tau_{q_0})$, $t \geq 0$. Using (5.35), we conclude that

$$
(5.37) \quad \begin{aligned}
&\iint F(V(\cdot))F^{(N)}(V^{(q_0)}(t_1),\ldots,V^{(q_0)}(t_N))\,d\mu \\
&= \iint F(\theta_{\tau_q}(V(\cdot)))F^{(N)}(V^{(q_0)}(t_1),\ldots,V^{(q_0)}(t_N))\,d\mu.
\end{aligned}
$$

Using Proposition 4.6 and an argument analogous to the one applied in the proof of Theorem 5.1 [see in particular the argument leading up to (5.4)], we conclude, upon the subsequent passages to the limit in (5.37), first as $q \to +\infty$, then as $q_0 \to +\infty$, that

$$
(5.38) \quad \begin{aligned}
&\iint F(V(\cdot))F^{(N)}(V(t_1),\ldots,V(t_N))\,d\mu \\
&= \iint F(V(\cdot))\,d\mu \times \iint F^{(N)}(V(t_1),\ldots,V(t_N))\,d\mu.
\end{aligned}
$$

From (5.36) and (5.38), we conclude that

$$
\left| \iint [F(V(\cdot))]^2 \, d\mu - \left[\iint F(V(\cdot))\,d\mu\right]^2 \right| < 2\varepsilon \sup |F|.
$$

Hence, upon the application of the fact that $\varepsilon > 0$ has been chosen arbitrarily, we get

$$
\iint [F(V(\cdot))]^2 \, d\mu = \left[\iint F(V(\cdot))\,d\mu\right]^2
$$

so $F(V(\cdot)) \equiv \text{const}$, $\mu$-a.s. $\square$

**6. The proof of the law of large numbers.** From (2.2) and Proposition 5.8, we immediately conclude that

$$
(6.1) \quad \lim_{t \to +\infty} \frac{\pi(t)}{t} = \lim_{t \to +\infty} \frac{1}{t} \int_0^t V(s)\,ds = \mathbf{v}_*, \qquad \mu\text{-a.s.},
$$

with $\mathbf{v}_*$ given by (1.2).

We show that the limit in (6.1) holds $P_0$-a.s. To demonstrate this fact, it suffices only to show that $\pi(t)/t$ converges $P_0$-a.s., as $t \to +\infty$, to a deterministic limit, which, as a consequence of the absolute continuity of $\mu$ w.r.t. $P_0$, must be equal to $\mathbf{v}_*$. In fact it suffices only to show the $P_0$-a.s. convergence of the sequence $(\pi(n)/n)_{n \geq 1}$. Indeed, we have

$$
(6.2) \quad \lim_{n \to +\infty} \sup_{n \leq t \leq n+1} \left| \frac{\pi(t)}{t} - \frac{\pi(n)}{n} \right| = 0.
$$

The latter is a consequence of the following estimate, that is an immediate consequence of (1.1):

$$(6.3) \quad P_0\left[\sup_{n\geq N}\sup_{n\leq t\leq n+1}\frac{|\pi(t)-\pi(n)|}{n}\geq\varepsilon\right]\leq \mathbb{W}\left[\sup_{t\geq N}\frac{|\pi(t)|}{[t]}\geq\frac{\varepsilon}{2}\right].$$

Inequality (6.3) holds for any $\varepsilon > 0$ and $N \geq 2(U+|\mathbf{v}|)/\varepsilon$, with $U$ the constant from condition (R). The right-hand side of (6.3) is as small as we wish, provided that $N$ is chosen sufficiently large.

Let $t_*$ be defined by (5.13). We also denote

$$\mathbf{w}_* := \int\int \pi(\tau_1)\,dP_{H_*}.$$

Note that, in consequence of Corollary 5.4 and Lemma 5.5 we have the following.

$$(6.4) \quad \lim_{n\to+\infty}\frac{\tau_n}{n}=t_*, \qquad \lim_{n\to+\infty}\frac{\pi(\tau_n)}{n}=\mathbf{w}_*, \qquad P_{H_*}\text{-a.s.}$$

Let us consider a nondecreasing sequence $(l_n)_{n\geq 1}$, that tends to infinity $P_{H_*}$-a.s., defined by

$$(6.5) \quad \tau_{l_n} \leq n < \tau_{l_n+1}.$$

We have

$$(6.6) \quad \lim_{n\uparrow+\infty}\frac{n}{l_n}=t_*, \qquad P_{H_*}\text{-a.s.}$$

Writing

$$(6.7) \quad \frac{\pi(n)}{n} = \frac{\pi(\tau_{l_n})}{l_n}\cdot\frac{l_n}{n} + \frac{\pi(n)-\pi(\tau_{l_n})}{n}$$

we conclude, by virtue of the individual ergodic theorem, that

$$(6.8) \quad \lim_{n\uparrow+\infty}\frac{\pi(n)}{n}=\frac{\mathbf{w}_*}{t_*}, \qquad P_{H_*}\text{-a.s.}$$

Let

$$E := \left[(\mathbf{x}_n,t_n)_{n\geq 1}\in(\mathbb{R}^d\times\mathbb{R})^{\mathbb{N}}:\frac{\sum_{m=1}^n \mathbf{x}_m}{\sum_{m=1}^n t_m}\not\to\frac{\mathbf{w}_*}{t_*},\right.$$

$$\left.\text{or }\frac{\sum_{m=1}^n t_m}{n}\not\to t_* \text{ as } n\uparrow+\infty\right].$$

From (6.4), it follows that

$$(6.9) \quad \int \mathbf{1}_E((\pi(\tau_{n+1})-\pi(\tau_n),\tau_{n+1}-\tau_n)_{n\geq 1})\,dP_{H_*}=0,$$

hence

$$\int \mathbf{M}_\omega[\mathbf{1}_E((\pi(\tau_{n+1}) - \pi(\tau_n), \tau_{n+1} - \tau_n)_{n \geq 1}), D = +\infty] H_*(\omega) \mathbb{P}(d\omega) = 0.$$

Since $H_* > 0$, $\mathbb{P}$-a.s., we conclude that

(6.10) $\mathbf{M}_\omega[\mathbf{1}_E((\pi(\tau_{n+1}) - \pi(\tau_n), \tau_{n+1} - \tau_n)_{n \geq 1}), D = +\infty] = 0,$ $\mathbb{P}$-a.s.

However, repeating the calculation made in (4.15) through (4.19), we obtain

(6.11)
$$\int \mathbf{1}_E((\pi(\tau_{n+2}) - \pi(\tau_{n+1}), \tau_{n+2} - \tau_{n+1})_{n \geq 1}) \, dP_0$$
$$= \int H(\omega') \mathbf{M}_{\omega'}[\mathbf{1}_E((\pi(\tau_{n+1}) - \pi(\tau_n), \tau_{n+1} - \tau_n)_{n \geq 1}),$$
$$D = +\infty] \mathbb{P}(d\omega')$$
$$\stackrel{(6.10)}{=} 0,$$

with

$$H(\omega') := \sum_{k,L=1}^{+\infty} \int_{\mathbb{R}^d} \int \mathcal{M}_{k,L,\mathbf{x}}(\omega, T_{-\mathbf{x}}\omega') \mathbb{F}_k(d\mathbf{x}) d\mathbb{P}(\omega).$$

From the definition of the set $E$, we conclude therefore

$$\lim_{n \uparrow +\infty} \frac{\pi(\tau_n)}{\tau_n} = \frac{\mathbf{w}_*}{t_*} \qquad \text{and} \qquad \lim_{n \uparrow +\infty} \frac{\tau_n}{n} = t_*, \qquad P_0\text{-a.s.}$$

Repeating the argument used in (6.5)–(6.6), this time with measure $P_0$ in place of $P_{H_*}$, we conclude that the limit in (6.8) holds $P_0$-a.s. $\square$

## APPENDIX A

**Proofs of Lemmas 3.1, 3.2 and Corollary 3.3.** Recall our standing assumption that $\kappa = 1$. Additionally, since all the estimates obtained below, as it becomes apparent in the course of the proofs, are independent of the choice of the starting point of the diffusion, we shall set $\mathbf{x} = \mathbf{0}$ throughout this section.

**A.1. Proof of (3.7).** For any $M > 0$, we denote

(A.1) $$\mathcal{S}_M^+ := [\mathbf{x} \in \mathbb{R}^d : -1 \leq \hat{\mathbf{v}} \cdot \mathbf{x} \leq M]$$

and $T_{\mathcal{S}_M^+}$ the exit time from the strip. Since

(A.2) $$Q_\omega[D = +\infty] = \lim_{M \uparrow +\infty} Q_\omega[T_{\mathcal{S}_M^+} < +\infty, \hat{\mathbf{v}} \cdot \pi(T_{\mathcal{S}_M^+}) = M],$$

inequality (3.7) will be proven once we show that there exists a constant $c > 0$, which bounds the right-hand side of (A.2). Recall that, for any $\omega \in \Omega$, the process $\mathbf{w}_\omega(\cdot)$, given by (2.2), is an $(\mathcal{M}_t)$ nonanticipative standard Brownian motion, under $Q_\omega$. Take any connected and bounded set $\mathbf{0} \in V \subseteq \mathcal{S}_M^+$. A simple argument involving the optional stopping theorem for martingales implies that

$$\text{(A.3)} \quad \mathbf{M}_\omega[\pi(n \wedge T_V) \cdot \hat{\mathbf{v}}] - \mathbf{M}_\omega\left[\int_0^{n \wedge T_V} \mathbf{u}(\pi(s)) \cdot \hat{\mathbf{v}} \, ds\right] = 0 \quad \forall n \geq 1.$$

Note that $-1 \leq \pi(n \wedge T_V) \cdot \hat{\mathbf{v}} \leq M$. On the other hand, with the choice of $\delta$ as in (2.1) we can write $\mathbf{M}_\omega T_V \leq (M+1)/\delta$ and in consequence $\mathbf{M}_\omega T_{\mathcal{S}_M^+} \leq (M+1)/\delta$, hence in particular

$$\text{(A.4)} \quad Q_\omega[T_{\mathcal{S}_M^+} < +\infty] = 1.$$

To finish the proof of (A.2), we choose $\theta_0 > 0$ so that $\theta_0 \leq \delta/2$ and let $\theta \in (0, \theta_0]$. Then, using the Markov property of $Q_\omega[\cdot]$ we get

$$\text{(A.5)} \quad \begin{aligned} &\mathbf{M}_\omega[\exp\{-\theta \hat{\mathbf{v}} \cdot \pi(t)\}|\mathcal{M}_s] \\ &= \exp\{-\theta \hat{\mathbf{v}} \cdot \pi(s)\} \\ &\quad + \int_s^t \mathbf{M}_\omega\{\exp\{-\theta \hat{\mathbf{v}} \cdot \pi(u)\}[-\theta \hat{\mathbf{v}} \cdot \mathbf{u}(\pi(u)) + \theta^2]|\mathcal{M}_s\} \, du \\ &\leq \exp\{-\theta \hat{\mathbf{v}} \cdot \pi(s)\} \end{aligned}$$

for $t \geq s$. The above calculation shows that $\exp\{-\theta \hat{\mathbf{v}} \cdot \pi(\cdot)\}$ is an $(\mathcal{M}_t)$-supermartingale. The optional sampling theorem for supermartingales and (A.4) yield

$$\mathbf{M}_\omega[\exp\{-\theta \hat{\mathbf{v}} \cdot \pi(\mathcal{S}_M^+)\}] \leq 1.$$

Thus, in consequence of the above estimate, we conclude that

$$e^\theta Q_\omega[\hat{\mathbf{v}} \cdot \pi(T_{\mathcal{S}_M^+}) = -1] \leq 1,$$

therefore

$$Q_\omega[\hat{\mathbf{v}} \cdot \pi(T_{\mathcal{S}_M^+}) = M] \geq 1 - e^{-\theta} \quad \forall M > 0,$$

and (3.7) is proven.

**A.2. Proof of (3.8).** We can write that the left hand-side of (3.8) is less than or equal to

$$\text{(A.6)} \quad Q_\omega[T_{\mathcal{S}_M} > t_M] + Q_\omega[T_{\mathcal{S}_M} \leq t_M, \hat{\mathbf{v}} \cdot \pi(T_{\mathcal{S}_M}) = -M].$$

Here $t_M := 2M\delta^{-1}$, $\mathcal{S}_M := [\mathbf{x} \in \mathbb{R}^d : -M \leq \hat{\mathbf{v}} \cdot \mathbf{x} \leq M]$ and $T_{\mathcal{S}_M}$ denotes the exit time from the strip $\mathcal{S}_M$. $\delta$ is defined in (2.1). Using the notation of (2.2), we can write that, on the event $[T_{\mathcal{S}_M} > t_M]$,

$$|\mathbf{w}_\omega(t_M)| = \left|\pi(t_M) - \int_0^{t_M} \mathbf{u}(\pi(s))\,ds\right| \geq M.$$

Hence,

$$Q_\omega[T_{\mathcal{S}_M} > t_M] \leq Q_\omega[|\mathbf{w}_\omega(t_M)| \geq M] \leq \exp\left\{-\frac{\delta M}{4}\right\}.$$

On the other hand,

(A.7) $\quad Q_\omega[T_{\mathcal{S}_M} \leq t_M, \hat{\mathbf{v}} \cdot \pi(T_{\mathcal{S}_M}) = -M] \leq Q_\omega\left[\sup_{0 \leq t \leq t_M} |\mathbf{w}_\omega(t)| \geq M\right].$

Using elementary estimates on the law of the maximum of a Brownian motion, we bound the right-hand side of (A.7) from above by $\exp\{-\frac{\delta M}{4d}\}$ and (3.8) follows.

**A.3. Proof of (3.9).** For any integer $m \geq 1$, we have

$$Q_\omega[2^m \leq M_* < 2^{m+1}, D < +\infty]$$

(A.8) $\quad \leq Q_\omega\left[|\pi(U_{2^m}) - 2^m \hat{\mathbf{v}}| \geq \frac{2^{m+1}|\mathbf{v}|}{\delta}\right]$

$$+ Q_\omega\left[|\pi(U_{2^m}) - 2^m \hat{\mathbf{v}}| < \frac{2^{m+1}|\mathbf{v}|}{\delta}, \tilde{U}_0 \circ \theta_{U_{2^m}} < U_{2^{m+1}} \circ \theta_{U_{2^m}}\right].$$

Let

$$C := \left[\mathbf{x} \in \mathbb{R}^d : |\mathbf{x} - (\hat{\mathbf{v}} \cdot \mathbf{x})\hat{\mathbf{v}}| \leq \frac{\mathbf{v} \cdot \mathbf{x}}{\delta}\right].$$

$C$ is a cone containing the support of the law of $\mathbf{u}(\mathbf{x})$, $\mathbf{x} \in \mathbb{R}^d$. Therefore $\int_0^t \mathbf{u}(\pi(s))\,ds \in C$, for all $t \geq 0$. On the other hand, there exists $c_1 > 0$ such that, for any $m \geq 1$, if

$$|\mathbf{x} - 2^m\hat{\mathbf{v}}| \geq \frac{2^{m+1}|\mathbf{v}|}{\delta} \quad \text{and} \quad \hat{\mathbf{v}} \cdot \mathbf{x} \leq 2^m,$$

then $\text{dist}(\mathbf{x}, C) > c_1 2^m$. The first term on the right-hand side of (A.8) can be therefore estimated by [since $\hat{\mathbf{v}} \cdot \pi(U_{2^m}) = 2^m$]

(A.9) $Q_\omega\left[|\pi(U_{2^m}) - 2^m\hat{\mathbf{v}}| \geq \frac{2^{m+1}|\mathbf{v}|}{\delta}, U_{2^m} \leq \frac{2^{m+1}}{\delta}\right] + Q_\omega\left[U_{2^m} > \frac{2^{m+1}}{\delta}\right].$

The expression in (A.9) can be therefore estimated by

$$Q_\omega\left[|\mathbf{w}_\omega(U_{2^m})| \geq c_1 2^m, U_{2^m} \leq \frac{2^{m+1}}{\delta}\right] + Q_\omega\left[\sup_{t \in [0, 2^{m+1}/\delta]} |\mathbf{w}_\omega(t)| \geq 2^m\right]$$

$$(A.10) \qquad \leq Q_\omega\left[\sup_{t \in [0, 2^{m+1}/\delta]} |\mathbf{w}_\omega(t)| \geq c_1 2^m, U_{2^m} \leq \frac{2^{m+1}}{\delta}\right]$$

$$+ Q_\omega\left[\sup_{t \in [0, 2^{m+1}/\delta]} |\mathbf{w}_\omega(t)| \geq 2^m\right].$$

Using once more the estimates on the law of the supremum of the Brownian motion, we bound the right-hand side of (A.10) from above by $c_2 \exp\{-c_3 2^m\}$ for some deterministic constants $c_2, c_3 > 0$ independent of $m$.

The second term on the right-hand side of (A.8) equals

$$\mathbf{M}_\omega\left[Q_{\pi(U_{2^m}),\omega}[\tilde{U}_0 < U_{2^{m+1}}], |\pi(U_{2^m}) - 2^m \hat{\mathbf{v}}| < \frac{2^{m+1}|\mathbf{v}|}{\delta}\right] \leq \exp\{-\gamma_1 2^m\},$$

by virtue of (the already proven) (3.8). We have therefore shown that

$$(A.11) \qquad Q_\omega[2^m \leq M_* < 2^{m+1}, D < +\infty] \leq c_4 \exp\{-c_5 2^m\}$$

for some deterministic constants $c_4, c_5 > 0$ independent of $m$, and (3.9) follows.

**A.4. Proof of Lemma 3.2.** For any $n \geq 1$, we obtain

$$(A.12) \qquad 0 = \mathbf{M}_\omega[\hat{\mathbf{v}} \cdot \mathbf{w}_\omega(n \wedge U_m)] \leq m - \delta \mathbf{M}_\omega(n \wedge U_m),$$

and (3.10) follows.

On the other hand,

$$Q_\omega[R_k < +\infty] \quad = \quad Q_\omega[S_k + D \circ \theta_{S_k} < +\infty]$$

$$\stackrel{\text{strong Markov prop.}}{=} \mathbf{M}_\omega[Q_{\pi(S_k),\omega}[D < +\infty], S_k < +\infty]$$

$$\stackrel{\text{Lemma 3.1}}{\leq} \quad (1-\gamma)Q_\omega[S_k < +\infty] \leq (1-\gamma)Q_\omega[R_{k-1} < +\infty]$$

and (3.11) follows by induction.

**A.5. Proof of Corollary 3.3.** Part (i) is an immediate conclusion from (3.11) and the Borel–Cantelli lemma. To show part (ii), note that

$$Q_{\mathbf{x},\omega}[S_K < +\infty]$$

$$= \sum_{k=1}^{+\infty} Q_{\mathbf{x},\omega}[R_{k-1} < +\infty, U_{M_{k-1}+r_0+1} \circ \theta_{R_{k-1}} < +\infty, K \circ \theta_{R_{k-1}} = 1]$$

(A.13)

$$= \sum_{k=1}^{+\infty} \int_{\mathbb{R}} \mathbf{M}_{\mathbf{x},\omega}[Q_{\pi(R_{k-1}),\omega}[U_{m+r_0+1} < +\infty, K = 1],$$

$$R_{k-1} < +\infty, M_{k-1} \in [m, m+dm)].$$

However, using (A.12), we can easily conclude that $Q_{\mathbf{y},\omega}[U_m < +\infty] = 1$ for all $\mathbf{y} \in \mathbb{R}^d$, $m \in \mathbb{R}$, $\mathbb{P}$-a.s.; hence the utmost right-hand side of (A.13) equals

$$\sum_{k=1}^{+\infty} \int_{\mathbb{R}} \mathbf{M}_{\mathbf{x},\omega}[Q_{\pi(R_{k-1}),\omega}[K=1], R_{k-1} < +\infty, M_{k-1} \in [m, m+dm)]$$

$$= Q_{\mathbf{x},\omega}[K < +\infty] = 1.$$

## APPENDIX B

**The existence of the isometric isomorphism $\mathcal{Z}$.** Suppose that $n \geq 1$ is a positive integer, $A_1, \ldots, A_n \in \mathcal{V}_0$, $B_1, \ldots, B_n \in \mathcal{R}$ are such that $A_1 \times B_1, \ldots, A_n \times B_n$ are pairwise disjoint and $c_1, \ldots, c_n \in \mathbb{R}$. We let

$$\mathcal{U}\left(\sum_{p=1}^{n} c_p \mathbf{1}_{A_p \times B_p}\right) := \sum_{p=1}^{n} c_p \mathbf{1}_{A_p} \mathbf{1}_{B_p}.$$

Since $\sigma$-algebras $\mathcal{V}_0$ and $\mathcal{R}$ are $\mathbb{P}$-independent, the mapping $\mathcal{U}$ is well defined and extends to a positivity-preserving isometry of any $L^p(\mathcal{T}_2 \otimes \mathcal{T}_3)$ into $L^p(\mathcal{T}_1)$ for any $p \in [1, +\infty]$. Thanks to the factorization property stated in the remark after condition (R), we conclude that $\mathcal{U}$ is in fact an isometric isomorphism between the relevant spaces. Define $\mathcal{Z} := \mathcal{U}^{-1}$. It is clear from the definition that properties (Z1) and (Z3) hold. Since $\mathcal{V}_t$ is generated by $\mathbb{P}$-independent $\sigma$-algebras $\mathcal{V}_0$ and $\mathcal{R}^t$, we can also immediately conclude (Z4). To prove condition (Z2,) we assume first that all $G_1, \ldots, G_N \in L^\infty(\mathcal{T}_2 \otimes \mathcal{T}_3)$. From the definition of $\mathcal{U}$, we conclude that $\mathcal{U}(G_1^{m_1} \cdots G_N^{m_N}) = [\mathcal{U}(G_1)]^{m_1} \cdots [\mathcal{U}(G_N)]^{m_N}$ for any nonnegative integers $m_1, \ldots, m_N \geq 0$. Hence, using, for example, the Weierstrass approximation theorem, we conclude that $\mathcal{U}\mathbf{\Phi}(G_1, \ldots, G_N) = \mathbf{\Phi}(\mathcal{U}(G_1), \ldots, \mathcal{U}(G_N))$ for any $\mathbf{\Phi} \in C_b(\mathbb{R}^N)$. We can remove the restriction on boundedness of $G_i$'s by using a standard truncation argument.

**Acknowledgment.** The authors would like to express their gratitude to the anonymous referee for a number of helpful comments that led to the improvement of the manuscript.

INVARIANT MEASURES FOR RANDOM DIFFUSIONS 33

INSTITUTE OF MATHEMATICS  
UMCS  
PL. MARII CURIE SKLODOWSKIEJ 1  
20-031 LUBLIN  
POLAND  
E-MAIL: komorow@golem.umcs.lublin.pl

FACULTY OF MATHEMATICS  
AND NATURAL SCIENCES  
CATHOLIC UNIVERSITY OF LUBLIN  
POLAND  
E-MAIL: grzegorz.krupa@kul.lublin.pl